\documentclass[leqno,12pt]{amsart}
\usepackage{amssymb} 
\theoremstyle{plain} 
\newtheorem{theorem}{\indent\sc Theorem}[section] 
\newtheorem{lemma}[theorem]{\indent\sc Lemma}

\newtheorem{proposition}[theorem]{\indent\sc Proposition}

\theoremstyle{definition} 
\newtheorem{definition}[theorem]{\indent\sc Definition}

\begin{document}

\title{Studies on the equations of Ince's table \\}
\author{Yusuke Sasano }

\renewcommand{\thefootnote}{\fnsymbol{footnote}}
\footnote[0]{2000\textit{ Mathematics Subjet Classification}.
34M55; 34M45; 58F05; 32S65.}

\maketitle

\begin{abstract}
We study the phase space of the equations of Ince's table from the viewpoint of its accessible singularities and local index.
\end{abstract}

\section{Introduction}

In 1979, K. Okamoto constructed the spaces of initial conditions of Painlev\'e equations, which can be considered as the parametrized spaces of all solutions, including the meromorphic solutions (see \cite{6}).

In this paper, we study the phase space of the equations of Ince's table (see \cite{Ince}) from the viewpoint of its accessible singularities and local index.

\section{Accessible singularity and local index}
Let us review the notion of {\it accessible singularity}. Let $B$ be a connected open domain in $\Bbb C$ and $\pi : {\mathcal W} \longrightarrow B$ a smooth proper holomorphic map. We assume that ${\mathcal H} \subset {\mathcal W}$ is a normal crossing divisor which is flat over $B$. Let us consider a rational vector field $\tilde v$ on $\mathcal W$ satisfying the condition
\begin{equation*}
\tilde v \in H^0({\mathcal W},\Theta_{\mathcal W}(-\log{\mathcal H})({\mathcal H})).
\end{equation*}
Fixing $t_0 \in B$ and $P \in {\mathcal W}_{t_0}$, we can take a local coordinate system $(x_1,\ldots ,x_n)$ of ${\mathcal W}_{t_0}$ centered at $P$ such that ${\mathcal H}_{\rm smooth \rm}$ can be defined by the local equation $x_1=0$.
Since $\tilde v \in H^0({\mathcal W},\Theta_{\mathcal W}(-\log{\mathcal H})({\mathcal H}))$, we can write down the vector field $\tilde v$ near $P=(0,\ldots ,0,t_0)$ as follows:
\begin{equation*}
\tilde v= \frac{\partial}{\partial t}+g_1 
\frac{\partial}{\partial x_1}+\frac{g_2}{x_1} 
\frac{\partial}{\partial x_2}+\cdots+\frac{g_n}{x_1} 
\frac{\partial}{\partial x_n}.
\end{equation*}
This vector field defines the following system of differential equations
\begin{equation}\label{39}
\frac{dx_1}{dt}=g_1(x_1,\ldots,x_n,t),\ \frac{dx_2}{dt}=\frac{g_2(x_1,\ldots,x_n,t)}{x_1},\cdots, \frac{dx_n}{dt}=\frac{g_n(x_1,\ldots,x_n,t)}{x_1}.
\end{equation}
Here $g_i(x_1,\ldots,x_n,t), \ i=1,2,\dots ,n,$ are holomorphic functions defined near $P=(0,\dots ,0,t_0).$

\begin{definition}\label{Def1}
With the above notation, assume that the rational vector field $\tilde v$ on $\mathcal W$ satisfies the condition
$$
(A) \quad \tilde v \in H^0({\mathcal W},\Theta_{\mathcal W}(-\log{\mathcal H})({\mathcal H})).
$$
We say that $\tilde v$ has an {\it accessible singularity} at $P=(0,\dots ,0,t_0)$ if
$$
x_1=0 \ {\rm and \rm} \ g_i(0,\ldots,0,t_0)=0 \ {\rm for \rm} \ {\rm every \rm} \ i, \ 2 \leq i \leq n.
$$
\end{definition}

If $P \in {\mathcal H}_{{\rm smooth \rm}}$ is not an accessible singularity, all solutions of the ordinary differential equation passing through $P$ are vertical solutions, that is, the solutions are contained in the fiber ${\mathcal W}_{t_0}$ over $t=t_0$. If $P \in {\mathcal H}_{\rm smooth \rm}$ is an accessible singularity, there may be a solution of \eqref{39} which passes through $P$ and goes into the interior ${\mathcal W}-{\mathcal H}$ of ${\mathcal W}$.

Here we review the notion of {\it local index}. Let $v$ be an algebraic vector field with an accessible singular point $\overrightarrow{p}=(0,\ldots,0)$ and $(x_1,\ldots,x_n)$ be a coordinate system in a neighborhood centered at $\overrightarrow{p}$. Assume that the system associated with $v$ near $\overrightarrow{p}$ can be written as
\begin{align}\label{b}
\begin{split}
\frac{d}{dt}\begin{pmatrix}
             x_1 \\
             x_2 \\
             \vdots\\
             x_{n-1} \\
             x_n
             \end{pmatrix}=\frac{1}{x_1}\left\{\begin{bmatrix}
             a_{11} & 0 & 0 & \hdots & 0 \\
             a_{21} & a_{22} & 0 &  \hdots & 0 \\
             \vdots & \vdots & \ddots & 0 & 0 \\
             a_{(n-1)1} & a_{(n-1)2} & \hdots & a_{(n-1)(n-1)} & 0 \\
             a_{n1} & a_{n2} & \hdots & a_{n(n-1)} & a_{nn}
             \end{bmatrix}\begin{pmatrix}
             x_1 \\
             x_2 \\
             \vdots\\
             x_{n-1} \\
             x_n
             \end{pmatrix}+\begin{pmatrix}
             x_1h_1(x_1,\ldots,x_n,t) \\
             h_2(x_1,\ldots,x_n,t) \\
             \vdots\\
             h_{n-1}(x_1,\ldots,x_n,t) \\
             h_n(x_1,\ldots,x_n,t)
             \end{pmatrix}\right\},\\
              (h_i \in {\Bbb C}(t)[x_1,\ldots,x_n], \ a_{ij} \in {\Bbb C}(t))
             \end{split}
             \end{align}
where $h_1$ is a polynomial which vanishes at $\overrightarrow{p}$ and $h_i$, $i=2,3,\ldots,n$ are polynomials of order at least 2 in $x_1,x_2,\ldots,x_n$, We call ordered set of the eigenvalues $(a_{11},a_{22},\cdots,a_{nn})$ {\it local index} at $\overrightarrow{p}$.

We are interested in the case with local index
\begin{equation}\label{integer}
(1,a_{22}/a_{11},\ldots,a_{nn}/a_{11}) \in {\Bbb Z}^{n}.
\end{equation}
These properties suggest the possibilities that $a_1$ is the residue of the formal Laurent series:
\begin{equation}
y_1(t)=\frac{a_{11}}{(t-t_0)}+b_1+b_2(t-t_0)+\cdots+b_n(t-t_0)^{n-1}+\cdots \quad (b_i \in {\Bbb C}),
\end{equation}
and the ratio $(1,a_{22}/a_{11},\ldots,a_{nn}/a_{11})$ is resonance data of the formal Laurent series of each $y_i(t) \ (i=2,\ldots,n)$, where $(y_1,\ldots,y_n)$ is original coordinate system satisfying $(x_1,\ldots,x_n)=(f_1(y_1,\ldots,y_n),\ldots,\\
f_n(y_1,\ldots,y_n)), \ f_i(y_1,\ldots,y_n) \in {\Bbb C}(t)(y_1,\ldots,y_n)$.

If each component of $(1,a_{22}/a_{11},\ldots,a_{nn}/a_{11})$ has the same sign, we may resolve the accessible singularity by blowing-up finitely many times. However, when different signs appear, we may need to both blow up and blow down.

The $\alpha$-test,
\begin{equation}\label{poiuy}
t=t_0+\alpha T, \quad x_i=\alpha X_i, \quad \alpha \rightarrow 0,
\end{equation}
yields the following reduced system:
\begin{align}\label{ppppppp}
\begin{split}
\frac{d}{dT}\begin{pmatrix}
             X_1 \\
             X_2 \\
             \vdots\\
             X_{n-1} \\
             X_n
             \end{pmatrix}=\frac{1}{X_1}\begin{bmatrix}
             a_{11}(t_0) & 0 & 0 & \hdots & 0 \\
             a_{21}(t_0) & a_{22}(t_0) & 0 &  \hdots & 0 \\
             \vdots & \vdots & \ddots & 0 & 0 \\
             a_{(n-1)1}(t_0) & a_{(n-1)2}(t_0) & \hdots & a_{(n-1)(n-1)}(t_0) & 0 \\
             a_{n1}(t_0) & a_{n2}(t_0) & \hdots & a_{n(n-1)}(t_0) & a_{nn}(t_0)
             \end{bmatrix}\begin{pmatrix}
             X_1 \\
             X_2 \\
             \vdots\\
             X_{n-1} \\
             X_n
             \end{pmatrix},
             \end{split}
             \end{align}
where $a_{ij}(t_0) \in {\Bbb C}$. Fixing $t=t_0$, this system is the system of the first order ordinary differential equation with constant coefficient. Let us solve this system. At first, we solve the first equation:
\begin{equation}
X_1(T)=a_{11}(t_0)T+C_1 \quad (C_1 \in {\Bbb C}).
\end{equation}
Substituting this into the second equation in \eqref{ppppppp}, we can obtain the first order linear ordinary differential equation:
\begin{equation}
\frac{dX_2}{dT}=\frac{a_{22}(t_0) X_2}{a_{11}(t_0)T+C_1}+a_{21}(t_0).
\end{equation}
By variation of constant, in the case of $a_{11}(t_0) \not= a_{22}(t_0)$ we can solve explicitly:
\begin{equation}
X_2(T)=C_2(a_{11}(t_0)T+C_1)^{\frac{a_{22}(t_0)}{a_{11}(t_0)}}+\frac{a_{21}(t_0)(a_{11}(t_0)T+C_1)}{a_{11}(t_0)-a_{22}(t_0)} \quad (C_2 \in {\Bbb C}).
\end{equation}
This solution is a single-valued solution if and only if
$$
\frac{a_{22}(t_0)}{a_{11}(t_0)} \in {\Bbb Z}.
$$
In the case of $a_{11}(t_0)=a_{22}(t_0)$ we can solve explicitly:
\begin{equation}
X_2(T)=C_2(a_{11}(t_0)T+C_1)+\frac{a_{21}(t_0)(a_{11}(t_0)T+C_1){\rm Log}(a_{11}(t_0)T+C_1)}{a_{11}(t_0)} \quad (C_2 \in {\Bbb C}).
\end{equation}
This solution is a single-valued solution if and only if
$$
a_{21}(t_0)=0.
$$
Of course, $\frac{a_{22}(t_0)}{a_{11}(t_0)}=1 \in {\Bbb Z}$.
In the same way, we can obtain the solutions for each variables $(X_3,\ldots,X_n)$. The conditions $\frac{a_{jj}(t)}{a_{11}(t)} \in {\Bbb Z}, \ (j=2,3,\ldots,n)$ are necessary condition in order to have the Painlev\'e property.

In the next section, in order to consider the phase spaces for each system, let us take the compactification $[z_0:z_1:z_2] \in {\Bbb P}^2$ of $(x,y) \in {\Bbb C}^2$ with the natural embedding
$$
(x,y)=(z_1/z_0,z_2/z_0).
$$
Moreover, we denote the boundary divisor in ${\Bbb P}^2$ by $ {\mathcal H}$. Extend the regular vector field on ${\Bbb C}^2$ to a rational vector field $\tilde v$ on ${\Bbb P}^2$. It is easy to see that ${\Bbb P}^2$ is covered by three copies of ${\Bbb C}^2${\rm : \rm}
\begin{align*}
&U_0={\Bbb C}^2 \ni (x,y),\\
&U_j={\Bbb C}^2 \ni (X_j,Y_j) \ (j=1,2),
\end{align*}
via the following rational transformations
\begin{align*}
& X_1=1/x, \quad Y_1=y/x,\\
& X_2=x/y, \quad Y_2=1/y,\\
\end{align*}

\section{Canonical equation of type I}

\begin{align}
\begin{split}
I:&\frac{d^2u}{dt^2}=0,\\
II:&\frac{d^2u}{dt^2}=6u^2, \quad III:\frac{d^2u}{dt^2}=6u^2+\frac{1}{2}, \quad IV:\frac{d^2u}{dt^2}=6u^2+t,\\
V:&\frac{d^2u}{dt^2}=-2u\frac{du}{dt}+q(t)\frac{du}{dt}+q'(t)u,\\
VI:&\frac{d^2u}{dt^2}=-3u\frac{du}{dt}-u^3+q(t)\left\{\frac{du}{dt}+u^2 \right\},\\
VII:&\frac{d^2u}{dt^2}=2u^3, \quad VIII:\frac{d^2u}{dt^2}=2u^3+\beta u+\gamma, \quad IX:\frac{d^2u}{dt^2}=2u^3+t u+\gamma,\\
X:&\frac{d^2u}{dt^2}=-u\frac{du}{dt}+u^3-12q(t) u+12q'(t),\\
\end{split}
\end{align}
where $':=\frac{d}{dt}$.

\section{Ince-V equation}
Ince-V equation is explicitly given by
\begin{equation}\label{eq;Ince-V}
\frac{d^2 u}{dt^2}=-2u\frac{du}{dt}+q(t)\frac{du}{dt}+q'(t)u.
\end{equation}
Here $u$ denotes unknown complex variable.

\begin{proposition}
The canonical transformation
\begin{equation*}
  \left\{
  \begin{aligned}
   x &=\frac{1}{u},\\
   y &=\frac{du}{dt}+u^2-q(t) u
   \end{aligned}
  \right. 
\end{equation*}
takes the equation \eqref{eq;Ince-V} to the system
\begin{equation}\label{system;Ince-V}
  \left\{
  \begin{aligned}
   \frac{dx}{dt} &=-x^2y-q(t)x+1,\\
   \frac{dy}{dt} &=0.
   \end{aligned}
  \right. 
\end{equation}
Here $x,y$ denote unknown complex variables.
\end{proposition}

\section{Ince-VI equation}
Ince-VI equation is explicitly given by
\begin{equation}\label{eq;Ince-VI}
\frac{d^2 u}{dt^2}=-3u\frac{du}{dt}-u^3+q(t)\left\{\frac{du}{dt}+u^2 \right\}.
\end{equation}

\begin{proposition}
The canonical transformation
\begin{equation*}
  \left\{
  \begin{aligned}
   x &=\frac{1}{u},\\
   y &=\frac{\frac{du}{dt}}{u}+u
   \end{aligned}
  \right. 
\end{equation*}
takes the equation \eqref{eq;Ince-VI} to the system
\begin{equation}\label{system;Ince-VI}
  \left\{
  \begin{aligned}
   \frac{dx}{dt} &=1-xy,\\
   \frac{dy}{dt} &=-y^2+q(t) y.
   \end{aligned}
  \right. 
\end{equation}
\end{proposition}
This system is a Riccati extension of the Riccati equation:
\begin{equation}
\frac{dy}{dt} =-y^2+q(t) y. 
\end{equation}

\section{Ince-VII equation}
Ince-VII equation is explicitly given by
\begin{equation}\label{eq;Ince-VII}
\frac{d^2 u}{dt^2}=2u^3.
\end{equation}

\begin{proposition}
The canonical transformation
\begin{equation*}
  \left\{
  \begin{aligned}
   x &=u,\\
   y &=\frac{\frac{du}{dt}}{u}
   \end{aligned}
  \right. 
\end{equation*}
takes the equation \eqref{eq;Ince-VII} to the system
\begin{equation}\label{system;Ince-VII}
  \left\{
  \begin{aligned}
   \frac{dx}{dt} &=xy,\\
   \frac{dy}{dt} &=2x^2-y^2.
   \end{aligned}
  \right. 
\end{equation}
\end{proposition}

\begin{figure}[h]
\unitlength 0.1in
\begin{picture}(44.00,22.40)(9.80,-23.00)
%
\special{pn 8}%
\special{pa 1180 610}%
\special{pa 2740 610}%
\special{fp}%
%
\special{pn 8}%
\special{pa 1410 180}%
\special{pa 2190 1820}%
\special{fp}%
%
\special{pn 8}%
\special{pa 2410 190}%
\special{pa 1820 1820}%
\special{fp}%
\put(9.8000,-12.8000){\makebox(0,0)[lb]{${\Bbb P}^2$}}%
%
\special{pn 20}%
\special{sh 0.600}%
\special{ar 2260 610 24 19  0.0000000 6.2831853}%
%
\special{pn 20}%
\special{sh 0.600}%
\special{ar 1910 600 24 19  0.0000000 6.2831853}%
%
\special{pn 20}%
\special{sh 0.600}%
\special{ar 1610 600 24 19  0.0000000 6.2831853}%
%
\special{pn 20}%
\special{pa 3610 590}%
\special{pa 5050 590}%
\special{fp}%
%
\special{pn 20}%
\special{pa 3950 290}%
\special{pa 3600 1050}%
\special{fp}%
\special{pa 3590 810}%
\special{pa 3890 1300}%
\special{fp}%
\special{pa 3900 1120}%
\special{pa 3440 1500}%
\special{fp}%
%
\special{pn 20}%
\special{pa 4400 260}%
\special{pa 4160 970}%
\special{fp}%
%
\special{pn 20}%
\special{pa 4810 280}%
\special{pa 5030 1030}%
\special{fp}%
%
\special{pn 20}%
\special{pa 5110 800}%
\special{pa 4750 1290}%
\special{fp}%
%
\special{pn 20}%
\special{pa 4730 1100}%
\special{pa 5250 1500}%
\special{fp}%
%
\special{pn 8}%
\special{pa 3560 1260}%
\special{pa 3430 1880}%
\special{dt 0.045}%
\special{pa 3430 1880}%
\special{pa 3430 1879}%
\special{dt 0.045}%
%
\special{pn 8}%
\special{pa 5090 1280}%
\special{pa 5330 1910}%
\special{dt 0.045}%
\special{pa 5330 1910}%
\special{pa 5330 1909}%
\special{dt 0.045}%
%
\special{pn 8}%
\special{pa 4170 770}%
\special{pa 4400 1270}%
\special{dt 0.045}%
\special{pa 4400 1270}%
\special{pa 4400 1269}%
\special{dt 0.045}%
%
\special{pn 8}%
\special{pa 3300 1700}%
\special{pa 4310 2300}%
\special{dt 0.045}%
\special{pa 4310 2300}%
\special{pa 4309 2300}%
\special{dt 0.045}%
\special{pa 4310 2300}%
\special{pa 5380 1690}%
\special{dt 0.045}%
\special{pa 5380 1690}%
\special{pa 5379 1690}%
\special{dt 0.045}%
%
\special{pn 20}%
\special{pa 3280 850}%
\special{pa 2820 850}%
\special{fp}%
\special{sh 1}%
\special{pa 2820 850}%
\special{pa 2887 870}%
\special{pa 2873 850}%
\special{pa 2887 830}%
\special{pa 2820 850}%
\special{fp}%
\put(36.1000,-2.3000){\makebox(0,0)[lb]{$E_7^{(1)}$-lattice}}%
\end{picture}%
\label{fig:Ince3}
\caption{Each bold line denotes $(-2)$-curve. The leftarrow denotes blowing-ups. The phase space $\mathcal X$ is the rational surface of type $E_7^{(1)}$.}
\end{figure}
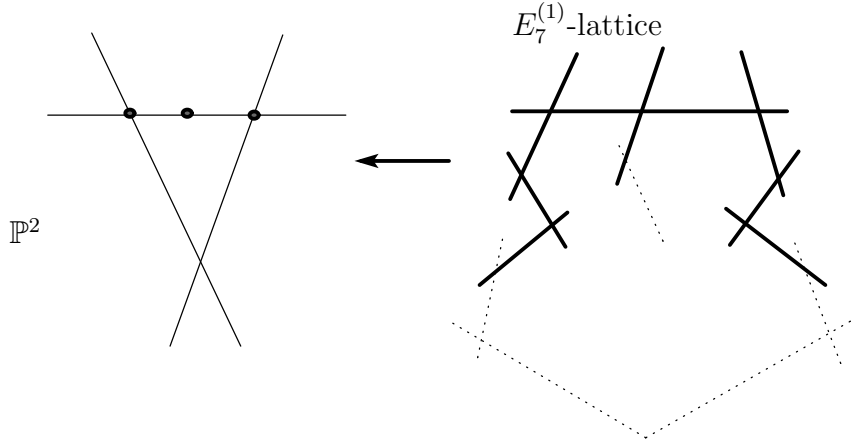

\begin{theorem}
After a series of explicit blowing-ups at ten points in ${\Bbb P}^2$, we obtain the rational surface ${\mathcal X}$ of type $E_7^{(1)}$ {\rm (see Figure 1) \rm}. The phase space ${\mathcal X}$ for the system \eqref{system;Ince-VII} is obtained by gluing four copies of ${\Bbb C}^2 \times {\Bbb C}${\rm:\rm}
\begin{center}
${U_j} \times {\Bbb C}={\Bbb C}^2 \times {\Bbb C} \ni \{(x_j,y_j,t)\},  \ \ j=0,1,2$
\end{center}
via the following birational transformations{\rm:\rm}
\begin{align*}
\begin{split}
0) \ &x_0=x, \quad y_0=y,\\
1) \ &x_1=xy, \quad y_1=\frac{1}{y},\\
2) \ &x_2=xy^3-y^4, \quad y_2=\frac{1}{y},\\
3) \ &x_3=xy^3+y^4, \quad y_3=\frac{1}{y}.
\end{split}
\end{align*}
\end{theorem}

\section{Ince-X equation}
Ince-X equation is explicitly given by
\begin{equation}\label{eq;Ince-X}
\frac{d^2 u}{dt^2}=-u\frac{du}{dt}+u^3-12q(t) u+12q'(t).
\end{equation}

\begin{proposition}
The canonical transformation
\begin{equation*}
  \left\{
  \begin{aligned}
   x &=u,\\
   y &=\frac{du}{dt}
   \end{aligned}
  \right. 
\end{equation*}
takes the equation \eqref{eq;Ince-X} to the system
\begin{equation}\label{system;Ince-X}
  \left\{
  \begin{aligned}
   \frac{dx}{dt} &=y,\\
   \frac{dy}{dt} &=-xy+x^3-12q(t) x+12q'(t).
   \end{aligned}
  \right. 
\end{equation}
\end{proposition}

\begin{theorem}\label{A}
After a series of explicit blowing-ups at eleven points in ${\Bbb P}^2$, we obtain the rational surface ${\mathcal X}$ of type $E_8^{(1)}$ {\rm (see Figure 2) \rm}. The phase space ${\mathcal X}$ for the system \eqref{system;Ince-X} is obtained by gluing three copies of ${\Bbb C}^2 \times {\Bbb C}${\rm:\rm}
\begin{center}
${U_j} \times {\Bbb C}={\Bbb C}^2 \times {\Bbb C} \ni \{(x_j,y_j,t)\},  \ \ j=0,1,2$
\end{center}
via the following birational transformations{\rm:\rm}
\begin{align*}
\begin{split}
0) \ &x_0=x, \quad y_0=y,\\
1) \ &x_1=\frac{1}{x}, \quad y_1=(y+x^2-12q(t))x,\\
2) \ &x_2=\frac{1}{x},\\
&y_2=\left\{(((y-x^2/2+6q(t))x-12q'(t))x+12(6q^2(t)-C_1 t-C_2)x-144q(t) q'(t) \right\}x.
\end{split}
\end{align*}
\end{theorem}

\begin{lemma}
The rational vector field $\tilde v$ has one accessible singular point{\rm : \rm}
\begin{equation*}
   P=\{(X_2,Y_2)|X_2=Y_2=0 \},
\end{equation*}
where $P$ has multiplicity of order $2$.
\end{lemma}

\begin{figure}[ht]
\unitlength 0.1in
\begin{picture}(53.30,42.30)(9.80,-44.10)
%
\special{pn 8}%
\special{pa 1180 610}%
\special{pa 2740 610}%
\special{fp}%
%
\special{pn 8}%
\special{pa 1410 180}%
\special{pa 2190 1820}%
\special{fp}%
%
\special{pn 8}%
\special{pa 2410 190}%
\special{pa 1820 1820}%
\special{fp}%
%
\special{pn 20}%
\special{pa 4670 810}%
\special{pa 4250 810}%
\special{fp}%
\special{sh 1}%
\special{pa 4250 810}%
\special{pa 4317 830}%
\special{pa 4303 810}%
\special{pa 4317 790}%
\special{pa 4250 810}%
\special{fp}%
%
\special{pn 8}%
\special{pa 4940 600}%
\special{pa 6070 600}%
\special{fp}%
%
\special{pn 8}%
\special{pa 5730 370}%
\special{pa 6300 1160}%
\special{fp}%
%
\special{pn 8}%
\special{pa 6310 960}%
\special{pa 6020 1690}%
\special{fp}%
%
\special{pn 20}%
\special{sh 0.600}%
\special{ar 6100 890 24 19  0.0000000 6.2831853}%
%
\special{pn 20}%
\special{sh 0.600}%
\special{ar 5970 700 24 19  0.0000000 6.2831853}%
%
\special{pn 20}%
\special{pa 1550 2440}%
\special{pa 1130 2440}%
\special{fp}%
\special{sh 1}%
\special{pa 1130 2440}%
\special{pa 1197 2460}%
\special{pa 1183 2440}%
\special{pa 1197 2420}%
\special{pa 1130 2440}%
\special{fp}%
%
\special{pn 8}%
\special{pa 1600 2240}%
\special{pa 2730 2240}%
\special{dt 0.045}%
\special{pa 2730 2240}%
\special{pa 2729 2240}%
\special{dt 0.045}%
%
\special{pn 8}%
\special{pa 2390 2010}%
\special{pa 2960 2800}%
\special{fp}%
%
\special{pn 20}%
\special{pa 2970 2600}%
\special{pa 2680 3330}%
\special{fp}%
%
\special{pn 20}%
\special{pa 2740 2330}%
\special{pa 2040 2570}%
\special{fp}%
%
\special{pn 20}%
\special{pa 2890 2450}%
\special{pa 2420 2950}%
\special{fp}%
\put(17.0000,-22.4000){\makebox(0,0)[lb]{$(-1)$}}%
\put(24.2000,-21.4000){\makebox(0,0)[lb]{$(-3)$}}%
\put(39.7000,-19.7000){\makebox(0,0)[lb]{$E_8^{(1)}$-lattice}}%
\put(9.8000,-12.8000){\makebox(0,0)[lb]{${\Bbb P}^2$}}%
%
\special{pn 20}%
\special{pa 3240 2470}%
\special{pa 3740 2470}%
\special{fp}%
\special{sh 1}%
\special{pa 3740 2470}%
\special{pa 3673 2450}%
\special{pa 3687 2470}%
\special{pa 3673 2490}%
\special{pa 3740 2470}%
\special{fp}%
%
\special{pn 8}%
\special{pa 3090 610}%
\special{pa 4110 610}%
\special{fp}%
%
\special{pn 8}%
\special{pa 3690 230}%
\special{pa 4160 1080}%
\special{fp}%
%
\special{pn 20}%
\special{sh 0.600}%
\special{ar 3900 610 24 19  0.0000000 6.2831853}%
%
\special{pn 20}%
\special{sh 0.600}%
\special{ar 2260 610 24 19  0.0000000 6.2831853}%
\put(46.0000,-40.0000){\makebox(0,0)[lb]{$\mathcal X$}}%
%
\special{pn 20}%
\special{pa 2170 2390}%
\special{pa 1970 2950}%
\special{fp}%
%
\special{pn 20}%
\special{pa 2500 2720}%
\special{pa 2420 3320}%
\special{fp}%
\special{pa 2510 3210}%
\special{pa 2070 3380}%
\special{fp}%
\special{pa 2130 3180}%
\special{pa 2130 3790}%
\special{fp}%
\special{pa 2200 3620}%
\special{pa 1770 3910}%
\special{fp}%
%
\special{pn 8}%
\special{pa 1550 2980}%
\special{pa 2150 2820}%
\special{dt 0.045}%
\special{pa 2150 2820}%
\special{pa 2149 2820}%
\special{dt 0.045}%
%
\special{pn 8}%
\special{pa 1850 3740}%
\special{pa 1930 4410}%
\special{dt 0.045}%
\special{pa 1930 4410}%
\special{pa 1930 4409}%
\special{dt 0.045}%
%
\special{pn 20}%
\special{pa 4770 1990}%
\special{pa 5340 2780}%
\special{fp}%
%
\special{pn 20}%
\special{pa 5350 2580}%
\special{pa 5060 3310}%
\special{fp}%
%
\special{pn 20}%
\special{pa 5120 2310}%
\special{pa 4420 2550}%
\special{fp}%
%
\special{pn 20}%
\special{pa 5270 2430}%
\special{pa 4800 2930}%
\special{fp}%
%
\special{pn 20}%
\special{pa 4550 2370}%
\special{pa 4350 2930}%
\special{fp}%
%
\special{pn 20}%
\special{pa 4880 2700}%
\special{pa 4800 3300}%
\special{fp}%
\special{pa 4890 3190}%
\special{pa 4450 3360}%
\special{fp}%
\special{pa 4510 3160}%
\special{pa 4510 3770}%
\special{fp}%
\special{pa 4580 3600}%
\special{pa 4150 3890}%
\special{fp}%
%
\special{pn 8}%
\special{pa 3930 2960}%
\special{pa 4530 2800}%
\special{dt 0.045}%
\special{pa 4530 2800}%
\special{pa 4529 2800}%
\special{dt 0.045}%
%
\special{pn 8}%
\special{pa 4230 3720}%
\special{pa 4310 4390}%
\special{dt 0.045}%
\special{pa 4310 4390}%
\special{pa 4310 4389}%
\special{dt 0.045}%
\end{picture}%
\label{fig:Ince2}
\caption{Each bold line denotes $(-2)$-curve. The leftarrow denotes blowing-ups and the rightarrow denotes blowing-downs. The symbol $(*)$ denotes intersection number of ${\Bbb P}^1$. The phase space $\mathcal X$ is the rational surface of type $E_8^{(1)}$.}
\end{figure}
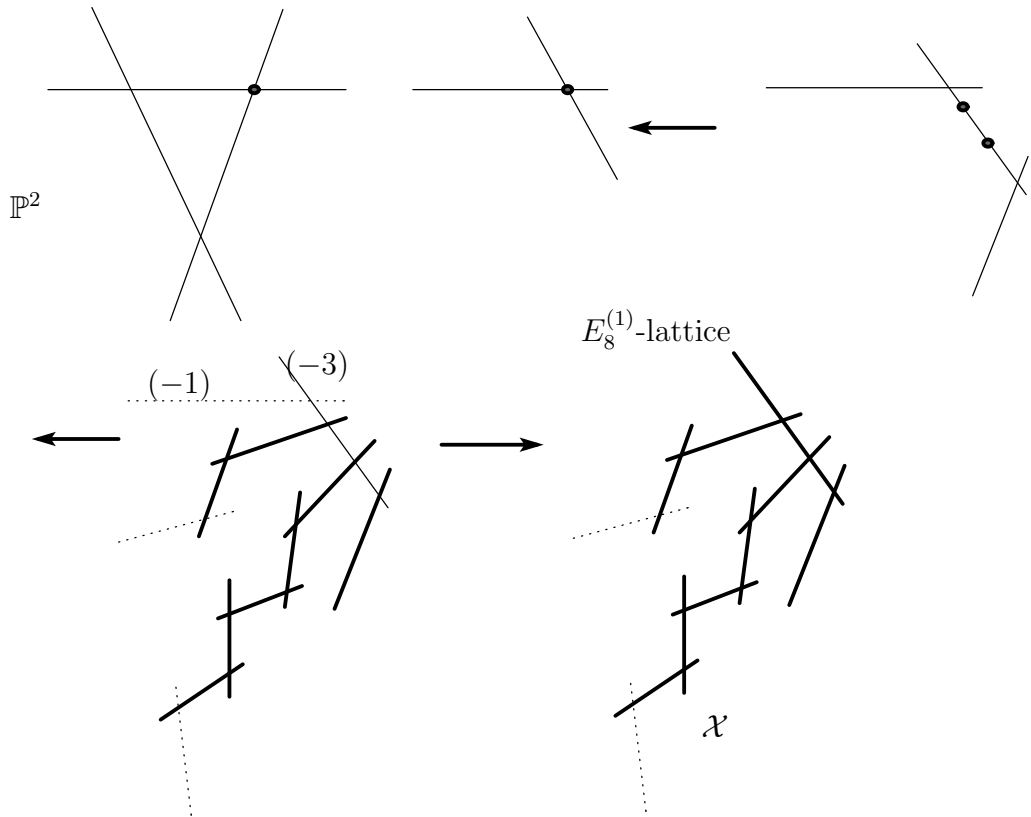

\section{Proof of Theorems \ref{A}}

Now we are ready to prove Theorem \ref{A}.

\vspace{0.5cm}
By the following steps, we can resolve the accessible singular point $P$.

{\bf Step 1}: We blow up at the point $P${\rm : \rm}
$$
x^{(1)}=X_2, \quad y^{(1)}=\frac{Y_2}{X_2}.
$$

{\bf Step 2}: We blow up at the point $\{(x^{(1)},y^{(1)})|x^{(1)}=y^{(1)}=0\}${\rm : \rm}
$$
x^{(2)}=\frac{x^{(1)}}{y^{(1)}}, \quad y^{(2)}=y^{(1)}.
$$

{\bf Step 3}: We make a chage of variables{\rm : \rm}
$$
x^{(3)}=\frac{1}{x^{(2)}}, \quad y^{(3)}=y^{(2)}.
$$

In the coordinate system $(x^{(3)},y^{(3)})$, the rational vector field $\tilde v$ has two accessible singular points{\rm : \rm}
\begin{align*}
\begin{split}
   P^{(1)} &=\{(x^{(3)},y^{(3)})|x^{(3)}=\frac{1}{2}, \ y^{(3)}=0\},\\
   P^{(2)} &=\{(x^{(3)},y^{(3)})|x^{(3)}=-1, \ y^{(3)}=0\}.
\end{split}
\end{align*}

Next let us calculate its local index at each point.
\begin{center}
\begin{tabular}{|c|c|c|} \hline 
Singular point & Type of local index & Condition \eqref{integer}   \\ \hline 
$P^{(1)}$ & $(-3,-\frac{1}{2})$ & $\frac{-3}{-\frac{1}{2}}=6$  \\ \hline 
$P^{(2)}$ & $(3,1)$ & $\frac{3}{1}=3$  \\ \hline
\end{tabular}
\end{center}
This property suggests that we will blow up six times to the direction $x^{(3)}$ on the resolution process of $P^{(1)}$ and three times to the direction $x^{(3)}$ on the resolution process of $P^{(2)}$.

At first, we will resolve the accessible singular point $P^{(1)}$.

{\bf Step 4}: We take the coordinate system centered at $P^{(1)}${\rm : \rm}
$$
x^{(4)}=x^{(3)}-\frac{1}{2}, \quad y^{(4)}=y^{(3)},
$$
and the system \eqref{system;Ince-X} is rewritten as follows:
\begin{align*}
\frac{d}{dt}\begin{pmatrix}
             x^{(4)} \\
             y^{(4)}
             \end{pmatrix}&=\frac{1}{y^{(4)}}\left\{\begin{pmatrix}
             -3 & 0 \\
             0 & -\frac{1}{2} 
             \end{pmatrix}\begin{pmatrix}
             x^{(4)} \\
             y^{(4)}
             \end{pmatrix}+\cdots\right\}
             \end{align*}
satisfying \eqref{b}.

{\bf Step 5}: We blow up at the point $\{(x^{(4)},y^{(4)})|x^{(4)}=y^{(4)}=0\}${\rm : \rm}
$$
x^{(5)}=\frac{x^{(4)}}{y^{(4)}}, \quad y^{(5)}=y^{(4)}.
$$

{\bf Step 6}: We blow up at the point $\{(x^{(5)},y^{(5)})|x^{(5)}=y^{(5)}=0\}${\rm : \rm}
$$
x^{(6)}=\frac{x^{(5)}}{y^{(5)}}, \quad y^{(6)}=y^{(5)}.
$$

{\bf Step 7}: We blow up at the point $\{(x^{(6)},y^{(6)})|x^{(6)}=-6q(t), \ y^{(6)}=0\}${\rm : \rm}
$$
x^{(7)}=\frac{x^{(6)}+6q(t)}{y^{(6)}}, \quad y^{(7)}=y^{(6)}.
$$

{\bf Step 8}: We blow up at the point $\{(x^{(7)},y^{(7)})|x^{(7)}=12q'(t), \ y^{(7)}=0\}${\rm : \rm}
$$
x^{(8)}=\frac{x^{(7)}-12q'(t)}{y^{(7)}}, \quad y^{(8)}=y^{(7)}.
$$

{\bf Step 9}: We blow up at the point $\{(x^{(8)},y^{(8)})|x^{(8)}=-12(6q^2(t)-C_1 t-C_2), \ y^{(8)}=0\}${\rm : \rm}
$$
x^{(9)}=\frac{x^{(8)}+12(6q^2(t)-C_1 t-C_2)}{y^{(8)}}, \quad y^{(9)}=y^{(8)}.
$$

{\bf Step 10}: We blow up at the point $\{(x^{(9)},y^{(9)})|x^{(9)}=144q(t) q'(t), \ y^{(9)}=0\}${\rm : \rm}
$$
x^{(10)}=\frac{x^{(9)}-144q(t) q'(t)}{y^{(9)}}, \quad y^{(10)}=y^{(9)},
$$
\begin{equation}\label{system1;Ince-X}
  \left\{
  \begin{aligned}
   \frac{dx^{(10)}}{dt} &=\frac{24\{12(q'(t))^2+12q(t)q''(t)-q''''(t) \}}{y^{(10)}}+g_1(x^{(10)},y^{(10)},t),\\
   \frac{dy^{(10)}}{dt} &=g_2(x^{(10)},y^{(10)},t),
   \end{aligned}
  \right. 
\end{equation}
where $g_i(x^{(10)},y^{(10)},t) \ (i=1,2)$ are polynomials in $x^{(10)},y^{(10)}$. Each right-hand side of the system \eqref{system1;Ince-X} is a {\it polynomial} if and only if
\begin{equation}
q''''(t)=12(q'(t))^2+12q(t)q''(t).
\end{equation}
This condition is reduced as follows:
\begin{equation}
q''(t)=6q^2(t)-C_1 t-C_2 \quad (C_1,C_2 \in {\Bbb C}).
\end{equation}
We remark that the coordinate system $(x^{(10)},y^{(10)})$ corresponds to $(y_2,x_2)$ given in Theorem \ref{A}.

Next, we will resolve the accessible singular point $P^{(2)}$.

{\bf Step 11}: We take the coordinate system centered at $P^{(2)}${\rm : \rm}
$$
x^{(11)}=x^{(3)}+1, \quad y^{(11)}=y^{(3)},
$$
and the system \eqref{system;Ince-X} is rewritten as follows:
\begin{align*}
\frac{d}{dt}\begin{pmatrix}
             x^{(11)} \\
             y^{(11)}
             \end{pmatrix}&=\frac{1}{y^{(11)}}\left\{\begin{pmatrix}
             3 & 0 \\
             0 & 1 
             \end{pmatrix}\begin{pmatrix}
             x^{(11)} \\
             y^{(11)}
             \end{pmatrix}+\cdots\right\}
             \end{align*}
satisfying \eqref{b}.

{\bf Step 12}: We blow up at the point $\{(x^{(11)},y^{(11)})|x^{(11)}=y^{(11)}=0\}${\rm : \rm}
$$
x^{(12)}=\frac{x^{(11)}}{y^{(11)}}, \quad y^{(12)}=y^{(11)}.
$$

{\bf Step 13}: We blow up at the point $\{(x^{(12)},y^{(12)})|x^{(12)}=12q(t), \ y^{(12)}=0\}${\rm : \rm}
$$
x^{(13)}=\frac{x^{(12)}-12q(t)}{y^{(12)}}, \quad y^{(13)}=y^{(12)}.
$$
Thus we have resolved the accessible singular point $P^{(2)}$. The coordinate system $(x^{(13)},y^{(13)})$ corresponds to $(y_1,x_1)$ given in Theorem \ref{A}.

Thus, we have completed the proof of Theorem \ref{A}.

\section{Canonical equation of type II}

\begin{align}
\begin{split}
XI:&\frac{d^2u}{dt^2}=\frac{1}{u} \left(\frac{du}{dt} \right)^2,\\
XII:&\frac{d^2u}{dt^2}=\frac{1}{u} \left(\frac{du}{dt} \right)^2+\alpha u^3+\beta u^2+\gamma +\frac{\delta}{u},\\
XIII:&\frac{d^2u}{dt^2}=\frac{1}{u} \left(\frac{du}{dt} \right)^2-\frac{1}{t}\frac{du}{dt}+\frac{1}{t}(\alpha u^2+\beta)+\gamma u^3+\frac{\delta}{u},\\
XIV:&\frac{d^2u}{dt^2}=\frac{1}{u} \left(\frac{du}{dt} \right)^2+\left\{q(t) u+\frac{r(t)}{u} \right\}+q'(t) u^2-r'(t),\\
XV:&\frac{d^2u}{dt^2}=\frac{1}{u} \left(\frac{du}{dt} \right)^2+\frac{1}{u}\frac{du}{dt}+r(t)u^2-u\frac{d}{dt} \left\{\frac{r'(t)}{r(t)} \right\},\\
XVI:&\frac{d^2u}{dt^2}=\frac{1}{u} \left(\frac{du}{dt} \right)^2-q'(t) \frac{1}{u}\frac{du}{dt}+u^3-q(t)u^2+q''(t),
\end{split}
\end{align}
where $':=\frac{d}{dt}$.

\section{Ince-XI equation}
Ince-XI equation is explicitly given by
\begin{equation}\label{eq;Ince-XI}
\frac{d^2 u}{dt^2}=\frac{1}{u}\left(\frac{du}{dt} \right)^2.
\end{equation}
Here $u$ denotes unknown complex variable.

\begin{proposition}
The canonical transformation
\begin{equation*}
  \left\{
  \begin{aligned}
   x &=u,\\
   y &=\frac{\frac{du}{dt}}{u}
   \end{aligned}
  \right. 
\end{equation*}
takes the equation \eqref{eq;Ince-XI} to the system
\begin{equation}\label{system;Ince-XI}
  \left\{
  \begin{aligned}
   \frac{dx}{dt} &=xy,\\
   \frac{dy}{dt} &=0.
   \end{aligned}
  \right. 
\end{equation}
Here $x,y$ denote unknown complex variables. This system can be solved by
\begin{equation}
(x,y)=(c_2e^{c_1 t},c_1) \quad (c_1,c_2 \in {\Bbb C}).
\end{equation}
\end{proposition}

\section{Ince-XII equation}
Ince-XII equation is explicitly given by
\begin{equation}\label{eq;Ince-XII}
\frac{d^2 u}{dt^2}=\frac{1}{u}\left(\frac{du}{dt} \right)^2+\alpha u^3+\beta u^2+\gamma+\frac{\delta}{u}.
\end{equation}

\begin{proposition}
The canonical transformation
\begin{equation*}
  \left\{
  \begin{aligned}
   x &=u,\\
   y &=\left(\frac{\frac{du}{dt}-\delta}{u}+\frac{\gamma}{\delta} \right)/u
   \end{aligned}
  \right. 
\end{equation*}
takes the equation \eqref{eq;Ince-XII} to the system
\begin{equation}\label{system;Ince-XII}
  \left\{
  \begin{aligned}
   \frac{dx}{dt} &=\frac{\partial H}{\partial y}=x^2y-\frac{\gamma}{\delta}x+\delta,\\
   \frac{dy}{dt} &=-\frac{\partial H}{\partial x}=-xy^2+\alpha x+\frac{\gamma}{\delta}y+\beta
   \end{aligned}
  \right. 
\end{equation}
with the polynomial Hamiltonian
\begin{equation}
H=\frac{x^2y^2}{2}-\frac{\alpha}{2}x^2-\frac{\gamma}{\delta}xy+\delta y-\beta x.
\end{equation}
\end{proposition}
We remark that the system \eqref{system;Ince-XII} has the Hamiltonian $H$ as its first integral.

\begin{theorem}
After a series of explicit blowing-ups in ${\Bbb P}^2$, we obtain the rational surface ${\mathcal X}$ of type $D_6^{(1)}$. The phase space ${\mathcal X}$ for the system \eqref{system;Ince-XII} can be obtained by gluing four copies of ${\Bbb C}^2 \times {\Bbb C}${\rm:\rm}
\begin{center}
${U_j} \times {\Bbb C}={\Bbb C}^2 \times {\Bbb C} \ni \{(x_j,y_j,t)\},  \ \ j=0,1,2,3$
\end{center}
via the following birational and symplectic transformations{\rm:\rm}
\begin{align*}
\begin{split}
0) \ &x_0=x, \quad y_0=y,\\
1) \ &x_1=\frac{1}{x}, \quad y_1=-\left((y-\sqrt{\alpha})x-\frac{\sqrt{\alpha} \gamma+\beta \delta}{\sqrt{\alpha} \delta} \right)x,\\
2) \ &x_2=\frac{1}{x}, \quad y_2=-\left((y+\sqrt{\alpha})x-\frac{\sqrt{\alpha} \gamma-\beta \delta}{\sqrt{\alpha} \delta} \right)x,\\
3) \ &x_3=x, \quad y_3=y-\frac{2\gamma}{\delta x}+\frac{2\delta}{x^2}.
\end{split}
\end{align*}
\end{theorem}

\begin{theorem}
The system \eqref{system;Ince-XII} admits the affine Weyl group symmetry of type $C_2^{(1)}$ as the group of its B{\"a}cklund transformations, whose generators $s_i$ are explicitly given as follows{\rm : \rm} with {\it the notation} $(*):=(x,y,t;\alpha,\beta,\gamma,\delta),$
\begin{align*}
\begin{split}
s_0:(x,y;\alpha,\beta,\gamma,\delta) \rightarrow & \left(x-\frac{\frac{\sqrt{\alpha} \gamma+\beta \delta}{\sqrt{\alpha} \delta}}{y-\sqrt{\alpha}},y;\alpha,-\frac{\gamma}{\delta} \sqrt{\alpha},-\frac{\beta \delta}{\sqrt{\alpha}},\delta \right),\\
s_1:(x,y;\alpha,\beta,\gamma,\delta) \rightarrow & \left(x-\frac{\frac{\sqrt{\alpha} \gamma-\beta \delta}{\sqrt{\alpha} \delta}}{y+\sqrt{\alpha}},y;\alpha,\frac{\gamma}{\delta} \sqrt{\alpha},\frac{\beta \delta}{\sqrt{\alpha}},\delta \right),\\
s_2:(x,y;\alpha,\beta,\gamma,\delta) \rightarrow & \left(x,y-\frac{2\gamma}{\delta x}+\frac{2\delta}{x^2};\alpha,\beta,\gamma,-\delta \right).
\end{split}
\end{align*}
\end{theorem}

\section{Ince-XIII equation}
Ince-XIII equation is explicitly given by
\begin{equation}\label{eq;Ince-XIII}
\frac{d^2 u}{dt^2}=\frac{1}{u} \left(\frac{du}{dt} \right)^2-\frac{1}{t}\frac{du}{dt}+\frac{1}{t}(\alpha u^2+\beta)+\gamma u^3+\frac{\delta}{u}.
\end{equation}

This equation coincides with the Painlev\'e III equation.

\section{Modified Ince-XIV equation}
Modified Ince-XIV equation is explicitly given by
\begin{equation}\label{eq;Ince-XIV}
\frac{d^2 u}{dt^2}=\frac{1}{u}\left(\frac{du}{dt} \right)^2+\frac{1}{2}\frac{d}{dt}\left(\frac{r'(t)}{r(t)} \right)u+\frac{r(t)}{u}+\frac{c_2 e^{c_1 t}}{\sqrt{r(t)}}u^2+c_1 \sqrt{r(t)}.
\end{equation}

\begin{proposition}
The canonical transformation
\begin{equation*}
  \left\{
  \begin{aligned}
   x &=u,\\
   y &=\left(\frac{\frac{du}{dt}-\sqrt{-r(t)}}{u}-\frac{2\sqrt{-1}c_1r(t)+r'(t)}{2r(t)} \right)/u
   \end{aligned}
  \right. 
\end{equation*}
takes the equation \eqref{eq;Ince-XIV} to the system
\begin{equation}\label{system;Ince-XIV}
  \left\{
  \begin{aligned}
   \frac{dx}{dt} &=\frac{\partial H}{\partial y}=x^2 y+\left(\sqrt{-1}c_1+\frac{r'(t)}{2r(t)} \right)x+\sqrt{-r(t)},\\
   \frac{dy}{dt} &=-\frac{\partial H}{\partial x}=-x y^2-\left(\sqrt{-1}c_1+\frac{r'(t)}{2r(t)} \right)y+\frac{c_2 e^{c_1 t}}{\sqrt{r(t)}}
   \end{aligned}
  \right. 
\end{equation}
with the polynomial Hamiltonian
\begin{equation}
H=\frac{x^2y^2}{2}+\left(\sqrt{-1}c_1+\frac{r'(t)}{2r(t)} \right)xy+\sqrt{-r(t)} y-\frac{c_2 e^{c_1 t}}{\sqrt{r(t)}}x.
\end{equation}
\end{proposition}

\begin{theorem}
After a series of explicit blowing-ups and blowing-downs in ${\Bbb P}^2$, we obtain the rational surface ${\mathcal X}$ of type $D_7^{(1)}$. The phase space ${\mathcal X}$ for the system \eqref{system;Ince-XIV} can be obtained by gluing three copies of ${\Bbb C}^2 \times {\Bbb C}${\rm:\rm}
\begin{center}
${U_j} \times {\Bbb C}={\Bbb C}^2 \times {\Bbb C} \ni \{(x_j,y_j,t)\},  \ \ j=0,1,2$
\end{center}
via the following birational and symplectic transformations{\rm:\rm}
\begin{align*}
\begin{split}
0) \ &x_0=x, \quad y_0=y,\\
1) \ &x_1=x+\frac{2(1+\sqrt{-1})c_1}{y}-\frac{2c_2 e^{c_1 t}}{\sqrt{r(t)} y^2}, \quad y_1=y,\\
2) \ &x_2=x, \quad y_2=y+\frac{2\sqrt{-1} c_1}{x}+\frac{2\sqrt{-r(t)}}{x^2}.
\end{split}
\end{align*}
\end{theorem}

\section{Ince-XV equation}
Ince-XV equation is explicitly given by
\begin{equation}\label{eq;Ince-XV}
\frac{d^2 u}{dt^2}=\frac{1}{u}\left( \frac{du}{dt} \right)^2-\frac{1}{u}\frac{du}{dt}+r(t) u^2-g(t) u.
\end{equation}

\begin{proposition}
The canonical transformation
\begin{equation*}
  \left\{
  \begin{aligned}
   x &=u,\\
   y &=\left(\frac{du}{dt}-1 \right)/u
   \end{aligned}
  \right. 
\end{equation*}
takes the equation \eqref{eq;Ince-XV} to the system
\begin{equation}\label{system;Ince-XV}
  \left\{
  \begin{aligned}
   \frac{dx}{dt} &=xy+1,\\
   \frac{dy}{dt} &=r(t) x-g(t),
   \end{aligned}
  \right. 
\end{equation}
where $r(t)$ and $g(t)$ satisfy the relation:
\begin{equation}
g(t)=\frac{r''(t)r(t)-(r'(t))^2}{r^2(t)}=\frac{d}{dt}\left\{\frac{r'(t)}{r(t)} \right\}.
\end{equation}
\end{proposition}

\begin{lemma}\label{lemma34}
The rational vector field $\tilde v$ has one accessible singular point{\rm : \rm}
\begin{equation*}
   P=\{(X_1,Y_1)|X_1=Y_1=0 \}.
\end{equation*}
\end{lemma}

By the following steps, we can resolve the accessible singular point $P$.

{\bf Step 1}: We blow up at the point $P${\rm : \rm}
$$
x^{(1)}=\frac{X_1}{Y_1}, \quad y^{(1)}=Y_1.
$$

{\bf Step 2}: We blow up at the point $\{(x^{(1)},y^{(1)})|x^{(1)}=y^{(1)}=0\}${\rm : \rm}
$$
x^{(2)}=x^{(1)}, \quad y^{(2)}=\frac{y^{(1)}}{x^{(1)}}.
$$

{\bf Step 3}: We make a chage of variables{\rm : \rm}
$$
x^{(3)}=x^{(2)}, \quad y^{(3)}=\frac{1}{y^{(2)}}.
$$

In the coordinate system $(x^{(3)},y^{(3)})$, the rational vector field $\tilde v$ has one accessible singular point{\rm : \rm}
\begin{align*}
\begin{split}
   P^{(1)} &=\{(x^{(3)},y^{(3)})|x^{(3)}=\frac{1}{2r(t)}, \ y^{(3)}=0 \}.
\end{split}
\end{align*}

Next let us calculate its local index at $P^{(1)}$.
\begin{center}
\begin{tabular}{|c|c|c|} \hline 
Singular point & Type of local index & Condition \eqref{integer}   \\ \hline 
$P^{(1)}$ & $(-1,-\frac{1}{2})$ & $\frac{-1}{-\frac{1}{2}}=2$  \\ \hline 
\end{tabular}
\end{center}
This property suggests that we will blow up second times to the direction $x^{(3)}$ on the resolution process of $P^{(1)}$.

{\bf Step 4}: We take the coordinate system centered at $P^{(1)}${\rm : \rm}
$$
x^{(4)}=x^{(3)}-\frac{1}{2r(t)}, \quad y^{(4)}=y^{(3)}.
$$

{\bf Step 5}: We blow up at the point $\{(x^{(4)},y^{(4)})|x^{(4)}=y^{(4)}=0\}${\rm : \rm}
$$
x^{(5)}=\frac{x^{(4)}}{y^{(4)}}, \quad y^{(5)}=y^{(4)}.
$$

{\bf Step 6}: We blow up at the point $\{(x^{(5)},y^{(5)})|x^{(5)}=\frac{r'(t)}{r^2(t)}, \ y^{(5)}=0\}${\rm : \rm}
$$
x^{(6)}=\frac{x^{(5)}-\frac{r'(t)}{r^2(t)}}{y^{(5)}}, \quad y^{(6)}=y^{(5)},
$$
\begin{equation}\label{system1;Ince-XV}
  \left\{
  \begin{aligned}
   \frac{dx^{(5)}}{dt} &=\left(g(t)-\frac{r''(t)r(t)-(r'(t))^2}{r^2(t)} \right)/y^{(5)}+g_1(x^{(5)},y^{(5)},t),\\
   \frac{dy^{(5)}}{dt} &=g_2(x^{(5)},y^{(5)},t),
   \end{aligned}
  \right. 
\end{equation}
where $g_i(x^{(5)},y^{(5)},t) \ (i=1,2)$ are polynomials in $x^{(5)},y^{(5)}$. Each right-hand side of the system \eqref{system1;Ince-XV} is a {\it polynomial} if and only if
\begin{equation}
g(t)=\frac{r''(t)r(t)-(r'(t))^2}{r^2(t)}=\frac{d}{dt}\left\{\frac{r'(t)}{r(t)} \right\}.
\end{equation}

Thus, we have resolved the accessible singular point $P^{(1)}$.

\section{Ince-XVI equation}
Ince-XVI equation is explicitly given by
\begin{equation}\label{eq;Ince-XVI}
\frac{d^2 u}{dt^2}=\frac{1}{u}\left(\frac{du}{dt} \right)^2-\frac{q'(t)}{u}\frac{du}{dt}+u^3-q(t) u^2+q''(t).
\end{equation}

\begin{proposition}
The canonical transformation
\begin{equation*}
  \left\{
  \begin{aligned}
   x &=u,\\
   y &=\frac{\frac{du}{dt}-q'(t)}{u}
   \end{aligned}
  \right. 
\end{equation*}
takes the equation \eqref{eq;Ince-XVI} to the system
\begin{equation}\label{system;Ince-XVI}
  \left\{
  \begin{aligned}
   \frac{dx}{dt} &=xy+q'(t),\\
   \frac{dy}{dt} &=x^2-q(t) x.
   \end{aligned}
  \right. 
\end{equation}
\end{proposition}

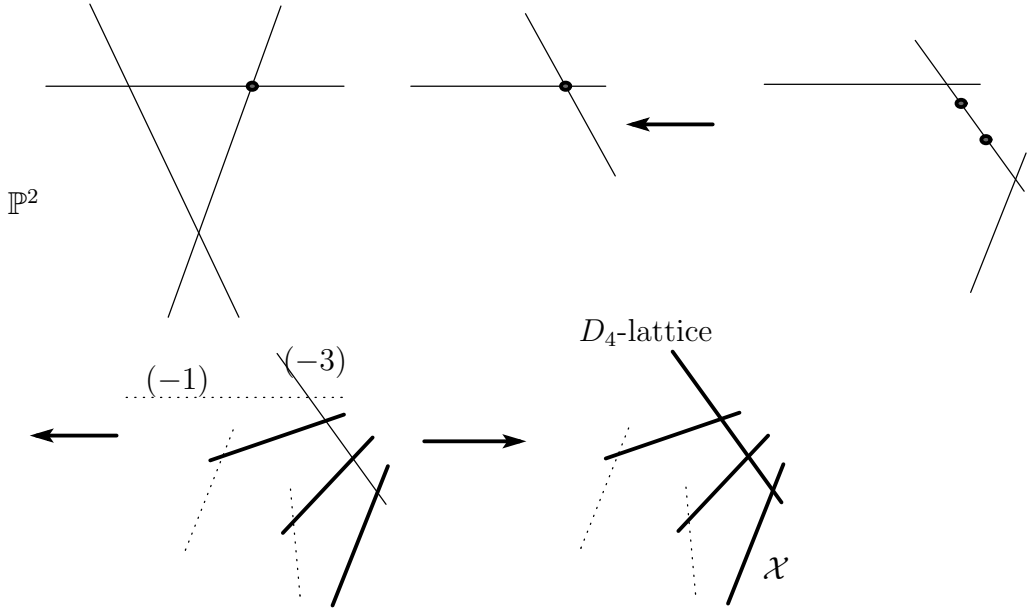
\begin{figure}
\unitlength 0.1in
\begin{picture}(53.30,31.50)(9.80,-33.30)
%
\special{pn 8}%
\special{pa 1180 610}%
\special{pa 2740 610}%
\special{fp}%
%
\special{pn 8}%
\special{pa 1410 180}%
\special{pa 2190 1820}%
\special{fp}%
%
\special{pn 8}%
\special{pa 2410 190}%
\special{pa 1820 1820}%
\special{fp}%
%
\special{pn 20}%
\special{pa 4670 810}%
\special{pa 4250 810}%
\special{fp}%
\special{sh 1}%
\special{pa 4250 810}%
\special{pa 4317 830}%
\special{pa 4303 810}%
\special{pa 4317 790}%
\special{pa 4250 810}%
\special{fp}%
%
\special{pn 8}%
\special{pa 4940 600}%
\special{pa 6070 600}%
\special{fp}%
%
\special{pn 8}%
\special{pa 5730 370}%
\special{pa 6300 1160}%
\special{fp}%
%
\special{pn 8}%
\special{pa 6310 960}%
\special{pa 6020 1690}%
\special{fp}%
%
\special{pn 20}%
\special{sh 0.600}%
\special{ar 6100 890 24 19  0.0000000 6.2831853}%
%
\special{pn 20}%
\special{sh 0.600}%
\special{ar 5970 700 24 19  0.0000000 6.2831853}%
%
\special{pn 20}%
\special{pa 1550 2440}%
\special{pa 1130 2440}%
\special{fp}%
\special{sh 1}%
\special{pa 1130 2440}%
\special{pa 1197 2460}%
\special{pa 1183 2440}%
\special{pa 1197 2420}%
\special{pa 1130 2440}%
\special{fp}%
%
\special{pn 8}%
\special{pa 1600 2240}%
\special{pa 2730 2240}%
\special{dt 0.045}%
\special{pa 2730 2240}%
\special{pa 2729 2240}%
\special{dt 0.045}%
%
\special{pn 8}%
\special{pa 2390 2010}%
\special{pa 2960 2800}%
\special{fp}%
%
\special{pn 20}%
\special{pa 2970 2600}%
\special{pa 2680 3330}%
\special{fp}%
%
\special{pn 20}%
\special{pa 2740 2330}%
\special{pa 2040 2570}%
\special{fp}%
%
\special{pn 20}%
\special{pa 2890 2450}%
\special{pa 2420 2950}%
\special{fp}%
%
\special{pn 8}%
\special{pa 2160 2410}%
\special{pa 1910 3040}%
\special{dt 0.045}%
\special{pa 1910 3040}%
\special{pa 1910 3039}%
\special{dt 0.045}%
%
\special{pn 8}%
\special{pa 2460 2720}%
\special{pa 2510 3280}%
\special{dt 0.045}%
\special{pa 2510 3280}%
\special{pa 2510 3279}%
\special{dt 0.045}%
\put(17.0000,-22.4000){\makebox(0,0)[lb]{$(-1)$}}%
\put(24.2000,-21.4000){\makebox(0,0)[lb]{$(-3)$}}%
%
\special{pn 20}%
\special{pa 4460 2000}%
\special{pa 5030 2790}%
\special{fp}%
%
\special{pn 20}%
\special{pa 5040 2590}%
\special{pa 4750 3320}%
\special{fp}%
%
\special{pn 20}%
\special{pa 4810 2320}%
\special{pa 4110 2560}%
\special{fp}%
%
\special{pn 20}%
\special{pa 4960 2440}%
\special{pa 4490 2940}%
\special{fp}%
%
\special{pn 8}%
\special{pa 4230 2400}%
\special{pa 3980 3030}%
\special{dt 0.045}%
\special{pa 3980 3030}%
\special{pa 3980 3029}%
\special{dt 0.045}%
%
\special{pn 8}%
\special{pa 4530 2710}%
\special{pa 4580 3270}%
\special{dt 0.045}%
\special{pa 4580 3270}%
\special{pa 4580 3269}%
\special{dt 0.045}%
\put(39.7000,-19.7000){\makebox(0,0)[lb]{$D_4$-lattice}}%
\put(9.8000,-12.8000){\makebox(0,0)[lb]{${\Bbb P}^2$}}%
%
\special{pn 20}%
\special{pa 3160 2470}%
\special{pa 3660 2470}%
\special{fp}%
\special{sh 1}%
\special{pa 3660 2470}%
\special{pa 3593 2450}%
\special{pa 3607 2470}%
\special{pa 3593 2490}%
\special{pa 3660 2470}%
\special{fp}%
%
\special{pn 8}%
\special{pa 3090 610}%
\special{pa 4110 610}%
\special{fp}%
%
\special{pn 8}%
\special{pa 3690 230}%
\special{pa 4160 1080}%
\special{fp}%
%
\special{pn 20}%
\special{sh 0.600}%
\special{ar 3900 610 24 19  0.0000000 6.2831853}%
%
\special{pn 20}%
\special{sh 0.600}%
\special{ar 2260 610 24 19  0.0000000 6.2831853}%
\put(49.3000,-31.9000){\makebox(0,0)[lb]{$\mathcal X$}}%
\end{picture}%
\label{fig:Ince1}
\caption{Each bold line denotes $(-2)$-curve. The leftarrow denotes blowing-ups and the rightarrow denotes blowing-downs. The symbol $(*)$ denotes intersection number of ${\Bbb P}^1$. The phase space $\mathcal X$ is the rational surface of type $D_4$.}
\end{figure}

\begin{theorem}
After a series of explicit blowing-ups and blowing-downs in ${\Bbb P}^2$, we obtain the rational surface ${\mathcal X}$ of type $D_4$  {\rm (see Figure 3) \rm}. The phase space ${\mathcal X}$ for the system \eqref{system;Ince-XVI} can be obtained by gluing three copies of ${\Bbb C}^2 \times {\Bbb C}${\rm:\rm}
\begin{center}
${U_j} \times {\Bbb C}={\Bbb C}^2 \times {\Bbb C} \ni \{(x_j,y_j,t)\},  \ \ j=0,1,2$
\end{center}
via the following birational transformations{\rm:\rm}
\begin{align*}
\begin{split}
0) \ &x_0=x, \quad y_0=y,\\
1) \ &x_1=\frac{1}{x}, \quad y_1=xy-x^2+q(t) x,\\
2) \ &x_2=\frac{1}{x}, \quad y_2=xy+x^2-q(t) x.
\end{split}
\end{align*}
\end{theorem}

\section{Canonical equation of type III}

\begin{align}
\begin{split}
XVII:&\frac{d^2u}{dt^2}=\frac{m-1}{mu} \left(\frac{du}{dt} \right)^2,\\
XVIII:&\frac{d^2u}{dt^2}=\frac{1}{2u} \left(\frac{du}{dt} \right)^2+4u^2,\\
XIX:&\frac{d^2u}{dt^2}=\frac{1}{2u} \left(\frac{du}{dt} \right)^2+4u^2+2u,\\
XX:&\frac{d^2u}{dt^2}=\frac{1}{2u} \left(\frac{du}{dt} \right)^2+4u^2+2tu,\\
XXI:&\frac{d^2u}{dt^2}=\frac{3}{4u} \left(\frac{du}{dt} \right)^2+3u^2,\\
XXII:&\frac{d^2u}{dt^2}=\frac{3}{4u} \left(\frac{du}{dt} \right)^2-1,\\
XXIII:&\frac{d^2u}{dt^2}=\frac{3}{4u} \left(\frac{du}{dt} \right)^2+3u^2+\alpha u+\beta,\\
XXIV:&\frac{d^2u}{dt^2}=\frac{m-1}{mu} \left(\frac{du}{dt} \right)^2+q(t)u\frac{du}{dt}-\frac{m q^2(t)}{(m+2)^2}u^3+\frac{m q'(t)}{m+2}u^2,\\
XXV:&\frac{d^2u}{dt^2}=\frac{3}{4u} \left(\frac{du}{dt} \right)^2-\frac{3u}{2}\frac{du}{dt}-\frac{u^3}{4}+\frac{q'(t)}{2q(t)} \left(u^2+\frac{du}{dt} \right)+r(t) u+q(t),\\
XXVI:&\frac{d^2u}{dt^2}=\frac{3}{4u} \left(\frac{du}{dt} \right)^2+\frac{6 q'(t)}{u}\frac{du}{dt}+3u^2+12q(t) u-12q''(t)-\frac{36 q'(t)^2}{u},\\
XXVII:&\frac{d^2u}{dt^2}=\frac{m-1}{mu} \left(\frac{du}{dt} \right)^2+\left(f(t) u+\phi(t)-\frac{m-2}{mu} \right)\frac{du}{dt}-\frac{m f(t)^2}{(m+2)^2}u^3\\
&+\frac{m(f'(t)-f(t) \phi(t))}{m+2}u^2+\psi(t) u-\phi(t)-\frac{1}{mu},\\
XXVIII:&\frac{d^2u}{dt^2}=\frac{1}{2u} \left(\frac{du}{dt} \right)^2-(u-q(t))\frac{du}{dt}+\frac{u^3}{2}-2q(t) u^2\\
&+3 \left(q'(t)+\frac{1}{2}q(t)^2 \right)u-\frac{72 r(t)^2}{u},\\
XXIX:&\frac{d^2u}{dt^2}=\frac{1}{2u} \left(\frac{du}{dt} \right)^2+\frac{3}{2} u^3,\\
XXX:&\frac{d^2u}{dt^2}=\frac{1}{2u} \left(\frac{du}{dt} \right)^2+\frac{3}{2} u^3+4\alpha u^2+2\beta u-\frac{\gamma^2}{2u},\\
XXXI:&\frac{d^2u}{dt^2}=\frac{1}{2u} \left(\frac{du}{dt} \right)^2+\frac{3}{2} u^3+4t u^2+2(t^2-\alpha) u-\frac{\beta^2}{2u},\\
XXXII:&\frac{d^2u}{dt^2}=\frac{1}{2u} \left(\frac{du}{dt} \right)^2-\frac{1}{2u},
\end{split}
\end{align}

\begin{align*}
\begin{split}
XXXIII:&\frac{d^2u}{dt^2}=\frac{1}{2u} \left(\frac{du}{dt} \right)^2+4u^2+\alpha u-\frac{1}{2u},\\
XXXIV:&\frac{d^2u}{dt^2}=\frac{1}{2u} \left(\frac{du}{dt} \right)^2+4\alpha u^2-t u-\frac{1}{2u},\\
XXXV:&\frac{d^2u}{dt^2}=\frac{2}{3u} \left(\frac{du}{dt} \right)^2-\left(\frac{2}{3}u-\frac{2}{3}q(t)-\frac{r(t)}{u} \right)\frac{du}{dt}+\frac{2}{3}u^3-\frac{10}{3}q(t) u^2\\
&+\left(4q'(t)+r(t)+\frac{8}{3}q(t)^2 \right)u-2q(t) r(t)-3r'(t)-\frac{3r(t)^2}{u},\\
XXXVI:&\frac{d^2u}{dt^2}=\frac{4}{5u} \left(\frac{du}{dt} \right)^2-\left(\frac{2}{5}u+\frac{1}{5}q(t)-\frac{r(t)}{u} \right)\frac{du}{dt}+\frac{4}{5}u^3+\frac{14}{5}q(t) u^2\\
&+\left(r(t)-3q'(t)+\frac{6}{5}q(t)^2 \right)u-\frac{1}{3}(q(t)r(t)+5r'(t))-\frac{5}{9}\frac{r(t)^2}{u},
\end{split}
\end{align*}
where $':=\frac{d}{dt}$.

\section{Ince-XVII equation}
Ince-XVII equation is explicitly given by
\begin{equation}\label{eq;Ince-XVII}
\frac{d^2 u}{dt^2}=\frac{m-1}{mu}\left(\frac{du}{dt} \right)^2.
\end{equation}
Here $u$ denotes unknown complex variable.

\begin{proposition}
The canonical transformation
\begin{equation*}
  \left\{
  \begin{aligned}
   x &=u,\\
   y &=\frac{\frac{du}{dt}}{u}
   \end{aligned}
  \right. 
\end{equation*}
takes the equation \eqref{eq;Ince-XVII} to the system
\begin{equation}\label{system;Ince-XVII}
  \left\{
  \begin{aligned}
   \frac{dx}{dt} &=xy,\\
   \frac{dy}{dt} &=-\frac{y^2}{m}.
   \end{aligned}
  \right. 
\end{equation}
Here $x,y$ denote unknown complex variables.
\end{proposition}

\section{Ince-XVIII equation}
Ince-XVIII equation is explicitly given by
\begin{equation}\label{eq;Ince-XVIII}
\frac{d^2 u}{dt^2}=\frac{1}{2u}\left(\frac{du}{dt} \right)^2+4u^2.
\end{equation}

\begin{proposition}
The canonical transformation
\begin{equation*}
  \left\{
  \begin{aligned}
   x &=u,\\
   y &=\frac{\frac{du}{dt}}{u}
   \end{aligned}
  \right. 
\end{equation*}
takes the equation \eqref{eq;Ince-XVIII} to the system
\begin{equation}\label{system;Ince-XVIII}
  \left\{
  \begin{aligned}
   \frac{dx}{dt} &=\frac{\partial H}{\partial y}=xy,\\
   \frac{dy}{dt} &=-\frac{\partial H}{\partial x}=-\frac{y^2}{2}+4x
   \end{aligned}
  \right. 
\end{equation}
with the polynomial Hamiltonian
\begin{equation*}
H=\frac{xy^2}{2}-2x^2.
\end{equation*}
\end{proposition}
We remark that the system \eqref{system;Ince-XVIII} has the Hamiltonian $H$ as its first integral.

\begin{theorem}
After a series of explicit blowing-ups in ${\Bbb P}^2$, we obtain the rational surface ${\mathcal X}$ of type $E_7^{(1)}$. The phase space ${\mathcal X}$ for the system \eqref{system;Ince-XVIII} can be obtained by gluing three copies of ${\Bbb C}^2 \times {\Bbb C}${\rm:\rm}
\begin{center}
${U_j} \times {\Bbb C}={\Bbb C}^2 \times {\Bbb C} \ni \{(x_j,y_j,t)\},  \ \ j=0,1,2$
\end{center}
via the following birational and symplectic transformations{\rm:\rm}
\begin{align*}
\begin{split}
0) \ &x_0=x, \quad y_0=y,\\
1) \ &x_1=-\left(x-\frac{y^2}{4} \right)y^2, \quad y_1=\frac{1}{y},\\
2) \ &x_2=-xy^2, \quad y_2=\frac{1}{y}.
\end{split}
\end{align*}
\end{theorem}

\section{Ince-XIX equation}
Ince-XIX equation is explicitly given by
\begin{equation}\label{eq;Ince-XIX}
\frac{d^2 u}{dt^2}=\frac{1}{2u}\left(\frac{du}{dt} \right)^2+4u^2+2u.
\end{equation}

\begin{proposition}
The canonical transformation
\begin{equation*}
  \left\{
  \begin{aligned}
   x &=u,\\
   y &=\frac{\frac{du}{dt}}{u}
   \end{aligned}
  \right. 
\end{equation*}
takes the equation \eqref{eq;Ince-XIX} to the system
\begin{equation}\label{system;Ince-XIX}
  \left\{
  \begin{aligned}
   \frac{dx}{dt} &=\frac{\partial H}{\partial y}=xy,\\
   \frac{dy}{dt} &=-\frac{\partial H}{\partial x}=-\frac{y^2}{2}+4x+2
   \end{aligned}
  \right. 
\end{equation}
with the polynomial Hamiltonian
\begin{equation*}
H=\frac{xy^2}{2}-2x^2-2x.
\end{equation*}
\end{proposition}
We remark that the system \eqref{system;Ince-XIX} has the Hamiltonian $H$ as its first integral.

\begin{theorem}
After a series of explicit blowing-ups in ${\Bbb P}^2$, we obtain the rational surface ${\mathcal X}$ of type $E_7^{(1)}$. The phase space ${\mathcal X}$ for the system \eqref{system;Ince-XIX} can be obtained by gluing three copies of ${\Bbb C}^2 \times {\Bbb C}${\rm:\rm}
\begin{center}
${U_j} \times {\Bbb C}={\Bbb C}^2 \times {\Bbb C} \ni \{(x_j,y_j,t)\},  \ \ j=0,1,2$
\end{center}
via the following birational and symplectic transformations{\rm:\rm}
\begin{align*}
\begin{split}
0) \ &x_0=x, \quad y_0=y,\\
1) \ &x_1=-\left(x-\frac{y^2}{4}+1 \right)y^2, \quad y_1=\frac{1}{y},\\
2) \ &x_2=-xy^2, \quad y_2=\frac{1}{y}.
\end{split}
\end{align*}
\end{theorem}

\section{Ince-XXI equation}
Ince-XXI equation is explicitly given by
\begin{equation}\label{eq;Ince-XXI}
\frac{d^2 u}{dt^2}=\frac{3}{4u}\left(\frac{du}{dt} \right)^2+3u^2.
\end{equation}

\begin{proposition}
The canonical transformation
\begin{equation*}
  \left\{
  \begin{aligned}
   x &=u,\\
   y &=\frac{\frac{du}{dt}}{u}
   \end{aligned}
  \right. 
\end{equation*}
takes the equation \eqref{eq;Ince-XXI} to the system
\begin{equation}\label{system;Ince-XXI}
  \left\{
  \begin{aligned}
   \frac{dx}{dt} &=xy,\\
   \frac{dy}{dt} &=-\frac{y^2}{4}+3x.
   \end{aligned}
  \right. 
\end{equation}
\end{proposition}
This system has its first integral:
\begin{equation}
I:=x(4x-y^2)^2.
\end{equation}

\begin{figure}[h]
\unitlength 0.1in
\begin{picture}(51.10,16.90)(9.80,-18.70)
%
\special{pn 8}%
\special{pa 1180 610}%
\special{pa 2740 610}%
\special{fp}%
%
\special{pn 8}%
\special{pa 1410 180}%
\special{pa 2190 1820}%
\special{fp}%
%
\special{pn 8}%
\special{pa 2410 190}%
\special{pa 1820 1820}%
\special{fp}%
\put(9.8000,-12.8000){\makebox(0,0)[lb]{${\Bbb P}^2$}}%
%
\special{pn 20}%
\special{sh 0.600}%
\special{ar 2260 610 24 19  0.0000000 6.2831853}%
%
\special{pn 20}%
\special{sh 0.600}%
\special{ar 1610 600 24 19  0.0000000 6.2831853}%
%
\special{pn 20}%
\special{pa 3280 850}%
\special{pa 2820 850}%
\special{fp}%
\special{sh 1}%
\special{pa 2820 850}%
\special{pa 2887 870}%
\special{pa 2873 850}%
\special{pa 2887 830}%
\special{pa 2820 850}%
\special{fp}%
%
\special{pn 20}%
\special{pa 3965 435}%
\special{pa 5699 435}%
\special{fp}%
%
\special{pn 20}%
\special{pa 4507 220}%
\special{pa 3760 849}%
\special{fp}%
%
\special{pn 20}%
\special{pa 5281 237}%
\special{pa 5716 689}%
\special{fp}%
\special{pa 5708 562}%
\special{pa 5343 871}%
\special{fp}%
\special{pa 5352 727}%
\special{pa 5743 1102}%
\special{fp}%
%
\special{pn 20}%
\special{pa 3902 567}%
\special{pa 4356 843}%
\special{fp}%
\special{pa 4302 738}%
\special{pa 4178 1113}%
\special{fp}%
%
\special{pn 8}%
\special{pa 4134 997}%
\special{pa 4658 1152}%
\special{dt 0.045}%
\special{pa 4658 1152}%
\special{pa 4657 1152}%
\special{dt 0.045}%
%
\special{pn 8}%
\special{pa 5530 1019}%
\special{pa 6090 1075}%
\special{dt 0.045}%
\special{pa 6090 1075}%
\special{pa 6089 1075}%
\special{dt 0.045}%
%
\special{pn 20}%
\special{pa 3790 630}%
\special{pa 3970 1400}%
\special{fp}%
%
\special{pn 8}%
\special{pa 3790 1200}%
\special{pa 4950 1870}%
\special{dt 0.045}%
\special{pa 4950 1870}%
\special{pa 4949 1870}%
\special{dt 0.045}%
\special{pa 4950 1870}%
\special{pa 6090 980}%
\special{dt 0.045}%
\special{pa 6090 980}%
\special{pa 6089 980}%
\special{dt 0.045}%
%
\special{pn 20}%
\special{pa 4950 1870}%
\special{pa 5240 1640}%
\special{fp}%
\special{sh 1}%
\special{pa 5240 1640}%
\special{pa 5175 1666}%
\special{pa 5198 1673}%
\special{pa 5200 1697}%
\special{pa 5240 1640}%
\special{fp}%
%
\special{pn 20}%
\special{pa 4960 1870}%
\special{pa 4590 1660}%
\special{fp}%
\special{sh 1}%
\special{pa 4590 1660}%
\special{pa 4638 1710}%
\special{pa 4636 1686}%
\special{pa 4658 1676}%
\special{pa 4590 1660}%
\special{fp}%
\end{picture}%
\label{fig:Ince4}
\caption{Each bold line denotes $(-2)$-curve. The leftarrow denotes blowing-ups.}
\end{figure}
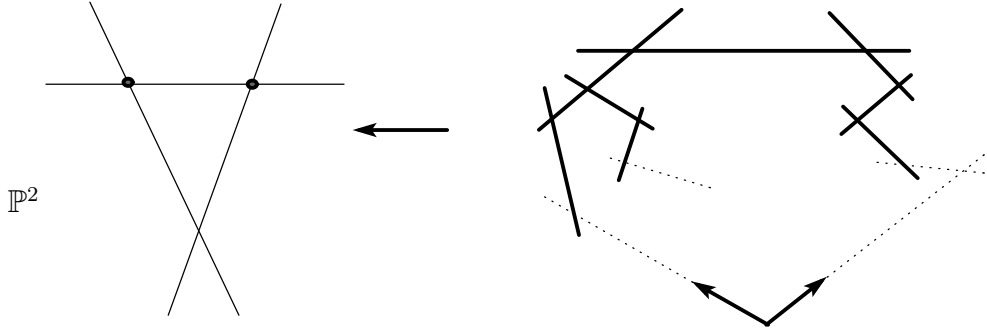

\begin{theorem}
The phase space ${\mathcal X}$ \rm{(see Figure 4) \rm} for the system \eqref{system;Ince-XXI} can be obtained by gluing three copies of ${\Bbb C}^2 \times {\Bbb C}${\rm:\rm}
\begin{center}
${U_j} \times {\Bbb C}={\Bbb C}^2 \times {\Bbb C} \ni \{(x_j,y_j,t)\},  \ \ j=0,1,2$
\end{center}
via the following birational transformations{\rm:\rm}
\begin{align*}
\begin{split}
0) \ &x_0=x, \quad y_0=y,\\
1) \ &x_1=\left(x-\frac{y^2}{4} \right)y, \quad y_1=\frac{1}{y},\\
2) \ &x_2=xy^4, \quad y_2=\frac{1}{y}.
\end{split}
\end{align*}
\end{theorem}

\section{Ince-XXII equation}
Ince-XXI equation is explicitly given by
\begin{equation}\label{eq;Ince-XXII}
\frac{d^2 u}{dt^2}=\frac{3}{4u}\left(\frac{du}{dt} \right)^2-1.
\end{equation}

\begin{proposition}
The canonical transformation
\begin{equation*}
  \left\{
  \begin{aligned}
   x &=\frac{u}{\frac{du}{dt}},\\
   y &=\frac{4u-\left(\frac{du}{dt} \right)^2}{4u\frac{du}{dt}}
   \end{aligned}
  \right. 
\end{equation*}
takes the equation \eqref{eq;Ince-XXII} to the system
\begin{equation}\label{system;Ince-XXII}
  \left\{
  \begin{aligned}
   \frac{dx}{dt} &=xy+\frac{1}{2},\\
   \frac{dy}{dt} &=y^2.
   \end{aligned}
  \right. 
\end{equation}
\end{proposition}
This system can be solved by
\begin{equation}
(x(t),y(t))=\left(\frac{t^2/4+c_1 t/2+c_2}{t+c_1},-\frac{1}{t+c_1} \right) \quad (c_1,c_2 \in {\Bbb C}).
\end{equation}

\section{Ince-XXIII equation}
Ince-XXIII equation is explicitly given by
\begin{equation}\label{eq;Ince-XXIII}
\frac{d^2 u}{dt^2}=\frac{1}{2u}\left(\frac{du}{dt} \right)^2+\frac{3}{2} u^3+\frac{\alpha}{2} u+\frac{\beta}{2u}.
\end{equation}

\begin{proposition}
The canonical transformation
\begin{equation*}
  \left\{
  \begin{aligned}
   x &=u,\\
   y &=\frac{\frac{du}{dt}-\sqrt{-\beta}}{u}
   \end{aligned}
  \right. 
\end{equation*}
takes the equation \eqref{eq;Ince-XXIII} to the system
\begin{equation}\label{system;Ince-XXIII}
  \left\{
  \begin{aligned}
   \frac{dx}{dt} &=\frac{\partial H}{\partial y}=xy+\sqrt{-\beta},\\
   \frac{dy}{dt} &=-\frac{\partial H}{\partial x}=\frac{3}{2}x^2-\frac{1}{2}y^2+\frac{\alpha}{2}
   \end{aligned}
  \right. 
\end{equation}
with the polynomial Hamiltonian
\begin{equation*}
H=\frac{1}{2}xy^2-\frac{1}{2}x^3-\frac{\alpha}{2}x+\sqrt{-\beta}y.
\end{equation*}
\end{proposition}
We remark that the system \eqref{system;Ince-XXIII} has the Hamiltonian $H$ as its first integral.

\begin{theorem}
After a series of explicit blowing-ups in ${\Bbb P}^2$, we obtain the rational surface ${\mathcal X}$ of type $E_6^{(1)}$. The phase space ${\mathcal X}$ for the system \eqref{system;Ince-XXIII} can be obtained by gluing four copies of ${\Bbb C}^2 \times {\Bbb C}${\rm:\rm}
\begin{center}
${U_j} \times {\Bbb C}={\Bbb C}^2 \times {\Bbb C} \ni \{(x_j,y_j,t)\},  \ \ j=0,1,2,3$
\end{center}
via the following birational and symplectic transformations{\rm:\rm}
\begin{align*}
\begin{split}
0) \ &x_0=x, \quad y_0=y,\\
1) \ &x_1=\frac{1}{x}, \quad y_1=-\left((y-x)x-\frac{\alpha-2\sqrt{-\beta}}{2} \right)x,\\
2) \ &x_2=\frac{1}{x}, \quad y_2=-\left((y+x)x+\frac{\alpha+2\sqrt{-\beta}}{2} \right)x,\\
3) \ &x_3=-(xy+2\sqrt{-\beta})y, \quad y_3=\frac{1}{y}.
\end{split}
\end{align*}
\end{theorem}

\begin{theorem}
The system \eqref{system;Ince-XXIII} is invariant under the following transformations\rm{: \rm}
\begin{align*}
\begin{split}
s_0:(x,y;\alpha,\beta) \rightarrow & \left(x-\frac{\frac{\alpha-2\sqrt{-\beta}}{2}}{y-x},y-\frac{\frac{\alpha-2\sqrt{-\beta}}{2}}{y-x};\frac{-\alpha+6\sqrt{-\beta}}{2},\frac{-\alpha^2-4\alpha \sqrt{-\beta}+4\beta}{16} \right),\\
s_1:(x,y;\alpha,\beta) \rightarrow & \left(x+\frac{\frac{\alpha+2\sqrt{-\beta}}{2}}{y+x},y-\frac{\frac{\alpha+2\sqrt{-\beta}}{2}}{y+x};\frac{-\alpha-6\sqrt{-\beta}}{2},\frac{-\alpha^2+4\alpha \sqrt{-\beta}+4\beta}{16} \right),\\
s_2:(x,y;\alpha,\beta) \rightarrow & \left(x,y+\frac{2\sqrt{-\beta}}{x};\alpha,\beta \right).
\end{split}
\end{align*}
\end{theorem}

\section{Ince-XXIV equation}
Ince-XXIV equation is explicitly given by
\begin{equation}\label{eq;Ince-XXIV}
\frac{d^2 u}{dt^2}=\frac{m-1}{mu}\left(\frac{du}{dt} \right)^2+q(t) u\frac{du}{dt}-\frac{mq^2(t)}{(m+2)^2} u^3+\frac{mq'(t)}{m+2} u^2.
\end{equation}

\begin{proposition}
The canonical transformation
\begin{equation*}
  \left\{
  \begin{aligned}
   x &=\frac{(m+1)(m+2)u}{-mq(t)u^2+(m+1)(m+2)\frac{du}{dt}},\\
   y &=\frac{(m+1)\left\{mq(t)u^2-(m+2)\frac{du}{dt} \right\}}{m^2q(t)u}
   \end{aligned}
  \right. 
\end{equation*}
takes the equation \eqref{eq;Ince-XXIV} to the system
\begin{equation}\label{system;Ince-XXIV}
  \left\{
  \begin{aligned}
   \frac{dx}{dt} &=-\frac{m^2q'(t)}{(m+1)(m+2)}x^2y-\frac{mq(t)}{m+1}xy-\frac{q'(t)}{q(t)}x-\frac{m+1}{m},\\
   \frac{dy}{dt} &=\frac{mq(t)}{(m+1)(m+2)}y^2-\frac{q'(t)}{q(t)}y.
   \end{aligned}
  \right. 
\end{equation}
\end{proposition}
This system is a Riccati extension of the Riccati equation:
\begin{equation}
\frac{dy}{dt}=\frac{mq(t)}{(m+1)(m+2)}y^2-\frac{q'(t)}{q(t)}y.
\end{equation}

\section{Ince-XXV equation}
Ince-XXV equation is explicitly given by
\begin{equation}\label{eq;Ince-XXV}
\frac{d^2 u}{dt^2}=\frac{3}{4u}\left(\frac{du}{dt} \right)^2-\frac{3}{2} u\frac{du}{dt}-\frac{1}{4} u^3+\frac{q'(t)}{2q(t)} \left(u^2+\frac{du}{dt} \right)+r(t)u+q(t).
\end{equation}

\begin{proposition}
The canonical transformation
\begin{equation*}
  \left\{
  \begin{aligned}
   x &=\frac{\frac{du}{dt}+u^2}{4q(t)u}+\frac{1}{\frac{du}{dt}+u^2},\\
   y &=\frac{u}{\frac{du}{dt}+u^2}
   \end{aligned}
  \right. 
\end{equation*}
takes the equation \eqref{eq;Ince-XXV} to the system
\begin{equation}\label{system;Ince-XXV}
  \left\{
  \begin{aligned}
   \frac{dx}{dt} &=-q(t) x^2-r(t) xy-\frac{q'(t)}{2q(t)}x+y+\frac{r(t)}{2q(t)},\\
   \frac{dy}{dt} &=-r(t) y^2-q(t) xy-\frac{q'(t)}{2q(t)}y+\frac{1}{2}.
   \end{aligned}
  \right. 
\end{equation}
\end{proposition}

\begin{proposition}
The phase space ${\mathcal X}$ for the system \eqref{system;Ince-XXV} is the $2$-dimensional projective space $(x,y) \in {\Bbb C}^2 \subset {\Bbb P}^2$.
\end{proposition}

\section{Ince-XXVI equation}
Ince-XXVI equation is explicitly given by
\begin{equation}\label{eq;Ince-XXVI}
\frac{d^2 u}{dt^2}=\frac{3}{4u}\left(\frac{du}{dt} \right)^2+\frac{6q'(t)}{u}\frac{du}{dt}+3u^2+12q(t) u-12q''(t)-\frac{36(q'(t))^2}{u}.
\end{equation}

\begin{proposition}
The canonical transformation
\begin{equation*}
  \left\{
  \begin{aligned}
   x &=u,\\
   y &=\frac{\frac{du}{dt}+12q'(t)}{u}
   \end{aligned}
  \right. 
\end{equation*}
takes the equation \eqref{eq;Ince-XXVI} to the system
\begin{equation}\label{system;Ince-XXVI}
  \left\{
  \begin{aligned}
   \frac{dx}{dt} &=xy-12q'(t),\\
   \frac{dy}{dt} &=-\frac{1}{4}y^2+3x+12q(t).
   \end{aligned}
  \right. 
\end{equation}
\end{proposition}

\begin{theorem}\label{th25}
After a series of explicit blowing-ups in ${\Bbb P}^2$, we obtain the rational surface ${\mathcal X}$ of type $E_8^{(1)}$. The phase space ${\mathcal X}$ for the system \eqref{system;Ince-XXVI} can be obtained by gluing three copies of ${\Bbb C}^2 \times {\Bbb C}${\rm:\rm}
\begin{center}
${U_j} \times {\Bbb C}={\Bbb C}^2 \times {\Bbb C} \ni \{(x_j,y_j,t)\},  \ \ j=0,1,2$
\end{center}
via the following birational transformations{\rm:\rm}
\begin{align*}
\begin{split}
0) \ &x_0=x, \quad y_0=y,\\
1) \ &x_1=\left(x-\frac{y^2}{4}+12q(t) \right)y, \quad y_1=\frac{1}{y},\\
2) \ &x_2=\left(((xy-16q'(t))y-32q''(t))y+768q(t)q'(t)-128q'''(t) \right)y, \quad y_2=\frac{1}{y}.
\end{split}
\end{align*}
\end{theorem}

\begin{lemma}
The rational vector field $\tilde v$ has two accessible singular points{\rm : \rm}
\begin{align*}
   P_1&=\{(X_1,Y_1)|X_1=Y_1=0 \},\\
   P_2&=\{(X_2,Y_2)|X_2=Y_2=0 \},
\end{align*}
where $P_1$ is multiple point of order $2$.
\end{lemma}

By the following steps, we can resolve the accessible singular point $P_2$.

Next let us calculate its local index at $P_2$.
\begin{center}
\begin{tabular}{|c|c|c|} \hline 
Singular point & Type of local index & Condition \eqref{integer}   \\ \hline 
$P_2$ & $(\frac{5}{4},\frac{1}{4})$ & $\frac{\frac{5}{4}}{\frac{1}{4}}=5$  \\ \hline 
\end{tabular}
\end{center}
This property suggests that we will blow up five times to the direction $X_1$ on the resolution process of $P_2$.

{\bf Step 1}: We blow up at the point $P_2${\rm : \rm}
$$
x^{(1)}=\frac{X_1}{Y_1}, \quad y^{(1)}=Y_1.
$$

{\bf Step 2}: We blow up at the point $\{(x^{(1)},y^{(1)})|x^{(1)}=y^{(1)}=0\}${\rm : \rm}
$$
x^{(2)}=\frac{x^{(1)}}{y^{(1)}}, \quad y^{(2)}=y^{(1)}.
$$

{\bf Step 3}: We blow up at the point $\{(x^{(2)},y^{(2)})|x^{(2)}=16q'(t), \ y^{(2)}=0\}${\rm : \rm}
$$
x^{(3)}=\frac{x^{(2)}-16q'(t)}{y^{(2)}}, \quad y^{(3)}=\frac{1}{y^{(2)}}.
$$

{\bf Step 4}: We blow up at the point $\{(x^{(3)},y^{(3)})|x^{(3)}=32q''(t), \ y^{(3)}=0\}${\rm : \rm}
$$
x^{(4)}=\frac{x^{(3)}-32q''(t)}{y^{(3)}}, \quad y^{(4)}=y^{(3)}.
$$

{\bf Step 5}: We blow up at the point $\{(x^{(4)},y^{(4)})|x^{(4)}=-768q(t)q'(t)+128q'''(t),\\
y^{(4)}=0\}${\rm : \rm}
$$
x^{(5)}=\frac{x^{(4)}+768q(t)q'(t)-128q'''(t)}{y^{(4)}}, \quad y^{(5)}=y^{(4)},
$$
\begin{equation}\label{4444}
  \left\{
  \begin{aligned}
   \frac{dx^{(5)}}{dt} &=\frac{128\left\{12((q'(t))^2+12q(t)q''(t)-q''''(t) \right\}}{y^{(5)}}+g_1(x^{(5)},y^{(5)},t),\\
   \frac{dy^{(5)}}{dt} &=g_2(x^{(5)},y^{(5)},t),
   \end{aligned}
  \right. 
\end{equation}
where $g_i(x^{(5)},y^{(5)},t) \ (i=1,2)$ are polynomials in $x^{(5)},y^{(5)}$. Each right-hand side of the system \eqref{4444} is a {\it polynomial} if and only if
\begin{equation}
q''''(t)=12((q'(t))^2+12q(t)q''(t).
\end{equation}
This equation reduces as follows:
\begin{equation}
q''(t)=6q^2(t)+C_1 t+C_2 \quad (C_1,C_2 \in {\Bbb C}).
\end{equation}

Thus, we have completed the proof of Theorem \ref{th25}.

\section{Ince-XXVII equation}
Ince-XXVII equation is explicitly given by
\begin{align}\label{eq;Ince-XXVII}
\begin {split}
\frac{d^2 u}{dt^2}=& \frac{m-1}{mu}\left(\frac{du}{dt} \right)^2+\left(fu+\phi-\frac{m-2}{mu} \right) \frac{du}{dt}\\
&-\frac{mf^2}{(m+2)^2}u^3+\frac{m(f'-f \phi)}{m+2}u^2+\psi u-\phi-\frac{1}{mu}.
\end{split}
\end{align}

\begin{proposition}
The canonical transformation
\begin{equation*}
  \left\{
  \begin{aligned}
   x &=\frac{1}{u},\\
   y &=\frac{\frac{du}{dt}-1}{u}-\frac{m}{m+2}f(t) u
   \end{aligned}
  \right. 
\end{equation*}
takes the equation \eqref{eq;Ince-XXVII} to the system
\begin{equation}\label{system;Ince-XXVII}
  \left\{
  \begin{aligned}
   \frac{dx}{dt} &=-x^2-xy-\frac{m}{m+2}f(t),\\
   \frac{dy}{dt} &=-\frac{y^2}{m}+\phi(t)y+f(t)+\psi(t)-\frac{m}{m+2}f(t).
   \end{aligned}
  \right. 
\end{equation}
\end{proposition}
This system is a Riccati extension of the Riccati equation:
\begin{equation}
\frac{dy}{dt}=-\frac{y^2}{m}+\phi(t)y+f(t)+\psi(t)-\frac{m}{m+2}f(t).
\end{equation}

\section{Ince-XXX equation}
Ince-XXX equation is explicitly given by
\begin{equation}\label{eq;Ince-XXX}
\frac{d^2 u}{dt^2}=\frac{1}{2u}\left(\frac{du}{dt} \right)^2+\frac{3}{2} u^3+4\alpha u^2+2\beta u-\frac{\gamma^2}{2u}.
\end{equation}

\begin{proposition}
The canonical transformation
\begin{equation*}
  \left\{
  \begin{aligned}
   x &=u,\\
   y &=\frac{\frac{du}{dt}-\gamma}{u}
   \end{aligned}
  \right. 
\end{equation*}
takes the equation \eqref{eq;Ince-XXX} to the system
\begin{equation}\label{system;Ince-XXX}
  \left\{
  \begin{aligned}
   \frac{dx}{dt} &=\frac{\partial H}{\partial y}=xy+\gamma,\\
   \frac{dy}{dt} &=-\frac{\partial H}{\partial x}=\frac{3}{2}x^2-\frac{1}{2}y^2+4\alpha x+2\beta
   \end{aligned}
  \right. 
\end{equation}
with the polynomial Hamiltonian
\begin{equation*}
H=\frac{1}{2}xy^2-\frac{1}{2}x^3-2\alpha x^2-2\beta x+\gamma y.
\end{equation*}
\end{proposition}
We remark that the system \eqref{system;Ince-XXX} has the Hamiltonian $H$ as its first integral.

\begin{theorem}
After a series of explicit blowing-ups in ${\Bbb P}^2$, we obtain the rational surface ${\mathcal X}$ of type $E_6^{(1)}$. The phase space ${\mathcal X}$ for the system \eqref{system;Ince-XXX} can be obtained by gluing four copies of ${\Bbb C}^2 \times {\Bbb C}${\rm:\rm}
\begin{center}
${U_j} \times {\Bbb C}={\Bbb C}^2 \times {\Bbb C} \ni \{(x_j,y_j,t)\},  \ \ j=0,1,2,3$
\end{center}
via the following birational and symplectic transformations{\rm:\rm}
\begin{align*}
\begin{split}
0) \ &x_0=x, \quad y_0=y,\\
1) \ &x_1=\frac{1}{x}, \quad y_1=-\left((y-x-2\alpha)x+2\alpha^2-2\beta+\gamma \right)x,\\
2) \ &x_2=\frac{1}{x}, \quad y_2=-\left((y+x+2\alpha)x-2\alpha^2+2\beta+\gamma \right)x,\\
3) \ &x_3=-(xy+2\gamma)y, \quad y_3=\frac{1}{y}.
\end{split}
\end{align*}
\end{theorem}

\begin{theorem}
The system \eqref{system;Ince-XXX} is invariant under the following transformations\rm{: \rm}
\begin{align*}
\begin{split}
s_0:(x,y;\alpha,\beta,\gamma) \rightarrow & \left(x+\frac{2\alpha^2-2\beta+\gamma}{y-x-2\alpha},y+\frac{2\alpha^2-2\beta+\gamma}{y-x-2\alpha};\alpha,\frac{2\alpha^2+\gamma}{2},\frac{-4\alpha^2+4\beta+\gamma}{3} \right),\\
s_1:(x,y;\alpha,\beta,\gamma) \rightarrow & (x+\frac{-2\alpha^2+2\beta+\gamma}{y+x+2\alpha},y-\frac{-2\alpha^2+2\beta+\gamma}{y+x+2\alpha};\\
&\alpha,\frac{6\alpha^2-2\beta-3\gamma}{4},\frac{2\alpha^2-2\beta+\gamma}{2}),\\
s_2:(x,y;\alpha,\beta,\gamma) \rightarrow & \left(x,y+\frac{2\gamma}{x};\alpha,\beta,-\gamma \right).
\end{split}
\end{align*}
\end{theorem}

\section{Ince-XXXII equation}
Ince-XXXII equation is explicitly given by
\begin{equation}\label{eq;Ince-XXXII}
\frac{d^2 u}{dt^2}=\frac{1}{2u}\left(\frac{du}{dt} \right)^2-\frac{1}{2u}.
\end{equation}

\begin{proposition}
The canonical transformation
\begin{equation*}
  \left\{
  \begin{aligned}
   x &=u,\\
   y &=\frac{\frac{du}{dt}-1}{u}
   \end{aligned}
  \right. 
\end{equation*}
takes the equation \eqref{eq;Ince-XXXII} to the system
\begin{equation}\label{system;Ince-XXXII}
  \left\{
  \begin{aligned}
   \frac{dx}{dt} &=xy+1,\\
   \frac{dy}{dt} &=-\frac{1}{2}y^2.
   \end{aligned}
  \right. 
\end{equation}
\end{proposition}
This system can be solved by
\begin{equation}
(x(t),y(t))=\left(-t+2c_1+c_2(t-2c_1)^2,\frac{2}{t-2c_1} \right) \quad (c_1,c_2 \in {\Bbb C}).
\end{equation}

\section{Ince-XXXV equation}
Ince-XXXV equation is explicitly given by
\begin{align}\label{eq;Ince-XXXV}
\begin{split}
\frac{d^2 u}{dt^2}=&\frac{2}{3u}\left(\frac{du}{dt} \right)^2-\left(\frac{2}{3}u-\frac{2}{3}q(t)-\frac{r(t)}{u} \right)\frac{du}{dt}+\frac{2}{3}u^3\\
&-\frac{10}{3}q(t)u^2+\left(4q'(t)+r(t)+\frac{8}{3}q^2(t) \right)u+2q(t)r(t)-3r'(t)-\frac{3r^2(t)}{u}.
\end{split}
\end{align}

\begin{proposition}
The canonical transformation
\begin{equation*}
  \left\{
  \begin{aligned}
   x &=u,\\
   y &=\frac{\frac{du}{dt}+3r(t)}{u}
   \end{aligned}
  \right. 
\end{equation*}
takes the equation \eqref{eq;Ince-XXXV} to the system
\begin{equation}\label{system;Ince-XXXV}
  \left\{
  \begin{aligned}
   \frac{dx}{dt} &=xy-3r(t),\\
   \frac{dy}{dt} &=\frac{2}{3}x^2-\frac{1}{3}y^2-\frac{2}{3}xy-\frac{10}{3}q(t)x+\frac{2}{3}q(t)y+\frac{8}{3}q^2(t)+3r(t)+4q'(t).
   \end{aligned}
  \right. 
\end{equation}
\end{proposition}

\begin{theorem}
After a series of explicit blowing-ups in ${\Bbb P}^2$, we obtain the rational surface ${\mathcal X}$ of type $E_7^{(1)}$. The phase space ${\mathcal X}$ for the system \eqref{system;Ince-XXXV} can be obtained by gluing four copies of ${\Bbb C}^2 \times {\Bbb C}${\rm:\rm}
\begin{center}
${U_j} \times {\Bbb C}={\Bbb C}^2 \times {\Bbb C} \ni \{(x_j,y_j,t)\},  \ \ j=0,1,2,3$
\end{center}
via the following birational transformations{\rm:\rm}
\begin{align*}
\begin{split}
0) \ &x_0=x, \quad y_0=y,\\
1) \ &x_1=\frac{1}{x}, \quad y_1=(y+x-4q(t))x,\\
2) \ &x_2=\frac{1}{x}, \quad y_2=-\left(((y-x/2+2q(t))x-9r(t)/2-6q'(t))x+9r'(t)+12q''(t) \right)x,\\
3) \ &x_3=-\left((xy-9q(t)/2)y+\frac{18q(t)r(t)-27r'(t)}{2} \right)y, \quad y_3=\frac{1}{y}.
\end{split}
\end{align*}
\end{theorem}

In the process of making its phase space, we can obtain the system in dimension three.

Setting
\begin{equation*}
x=q(t), \quad y=q'(t), \quad z=r(t),
\end{equation*}
we obtain the system
\begin{equation}\label{3rd}
  \left\{
  \begin{aligned}
   \frac{dx}{dt} &=y,\\
   \frac{dy}{dt} &=-2xy-3xz-\frac{3c_1+c_2}{4},\\
   \frac{dz}{dt} &=2xz+c_1.
   \end{aligned}
  \right. 
\end{equation}
Here $x,y,z$ denote unknown complex variables and $c_1,c_2$ are constant parameters.

\begin{theorem}
The system \eqref{3rd} is invariant under the following transformations\rm{: \rm}
\begin{align*}
\begin{split}
s_0:(x,y,z;c_1,c_2) \rightarrow & \left(x+\frac{c_1}{z},y-\frac{2c_1x}{z}-\frac{c_1^2}{z^2},z;-c_1,-6c_1+c_2 \right),\\
s_1:(x,y,z;c_1,c_2) \rightarrow & (\frac{4xy+3xz+c_2}{4y+3z}\\
&\frac{16y^3+24y^2z+9yz^2+8c_2xy+6c_2xz+c_2^2}{(4y+3z)^2},\\
&\frac{27z^3+48y^2z+72yz^2-32c_2xy-24c_2xz-4c_2^2}{3(4y+3z)^2};c_1-2c_2/3,-c_2).
\end{split}
\end{align*}
\end{theorem}
Here, the transformations $s_0,s_1$ satisfy the relations $s_0^2=s_1^2=1$.

\begin{theorem}
The phase space ${\mathcal X}$ for the system \eqref{3rd} can be obtained by gluing three copies of ${\Bbb C}^3 \times {\Bbb C}${\rm:\rm}
\begin{center}
${U_j} \times {\Bbb C}={\Bbb C}^3 \times {\Bbb C} \ni \{(x_j,y_j,z_j,t)\},  \ \ j=0,1,2$
\end{center}
via the following birational transformations{\rm:\rm}
\begin{align*}
\begin{split}
0) \ &x_0=x, \quad y_0=y, \quad z_0=z,\\
1) \ &x_1=\frac{1}{x}, \quad y_1=y+x^2, \quad z_1=-(xz+c_1)x,\\
2) \ &x_2=\frac{1}{x}, \quad y_2=\frac{1}{3}x(4xy+3xz+c_2), \quad z_2=\frac{2}{3}(2x^2+2y+3z).
\end{split}
\end{align*}
\end{theorem}

In the coordinate system $(x_2,y_2,z_2)$, we obtain the system
\begin{equation}\label{3rd1}
  \left\{
  \begin{aligned}
   \frac{dx_2}{dt} &=\frac{3}{2}x_2^4y_2+\frac{c_2}{2}x_2^3+\frac{3}{4}\left\{(c_1-c_2/3)t+c_3 \right\}x_2^2-1,\\
   \frac{dy_2}{dt} &=-3x_2^3y_2^2-\frac{3c_2}{2}x_2^2y_2-\frac{3}{2}\left\{(c_1-c_2/3)t+c_3 \right\}x_2y_2-\frac{c_2^2x_2}{6}-\frac{c_2}{4}\left\{(c_1-c_2/3)t+c_3 \right\},\\
   \frac{dz_2}{dt} &=c_1-\frac{1}{3}c_2.
   \end{aligned}
  \right. 
\end{equation}
We see that this system has its first integral:
\begin{equation}
\frac{dz_2}{dt}=c_1-\frac{1}{3}c_2,
\end{equation}
and we can solve explicitly as follows:
\begin{equation}
z_2(t)=\left(c_1-\frac{1}{3}c_2 \right)t+c_3,
\end{equation}
where $c_3$ is its integral constant.

Setting
\begin{equation*}
X=\frac{1}{x_2}, \quad Y=(x_2y_2+c_2/3)x_2,
\end{equation*}
we obtain the system
\begin{equation}
  \left\{
  \begin{aligned}
   \frac{dX}{dt} &=\frac{\partial H}{\partial Y}=X^2-\frac{3}{2}Y-\frac{1}{4}(3c_1-c_2)t-\frac{3}{4}c_3,\\
   \frac{dY}{dt} &=-\frac{\partial H}{\partial X}=-2XY+\frac{c_2}{3}
   \end{aligned}
  \right. 
\end{equation}
with the polynomial Hamiltonian
\begin{equation*}
H=X^2Y-\frac{3}{4}Y^2-\frac{1}{4}(3c_1-c_2)tY-\frac{3}{4}c_3Y-\frac{c_2}{3}X.
\end{equation*}
This system is the second Painlev\'e system.

\end{document}